\documentclass[]{elsarticle}

\usepackage[title]{appendix}
\usepackage{amsmath}
\usepackage{mathtools}
\usepackage{amssymb}
\usepackage{upgreek}
\usepackage{fullpage}
\usepackage{color}
\usepackage{caption}
\usepackage{subcaption}
\usepackage{makecell}
\usepackage{mathrsfs}
\usepackage{breqn}
\usepackage{float}
\usepackage{tikz}
\usetikzlibrary{arrows,shapes,backgrounds,decorations.markings}

\usepackage{titlesec}

\usepackage[linesnumbered,ruled,vlined]{algorithm2e}
\SetKwInput{KwInput}{Input}                
\SetKwInput{KwOutput}{Output}              

\usepackage{float}
\usepackage{booktabs}
\usepackage{tabu} 
\usepackage{dcolumn}
\newcolumntype{A}{D{.}{.}{2.3}}
\usepackage{tabularx}
\usepackage[hidelinks]{hyperref}

\usepackage{amsthm}

\newcommand{\bb}{\mathrm{b}}
\newcommand{\comp}{\mathrm{comp}}
\newcommand{\fin}{\mathrm{fin}}
\newcommand{\port}{\mathrm{port}}
\newcommand{\train}{\mathrm{train}}
\newcommand{\test}{\mathrm{test}}
\newcommand{\rb}{\mathrm{rb}}
\newcommand{\sample}{\mathrm{sample}}

\newcommand{\boplus}{\textstyle\bigoplus}

\DeclareMathOperator*{\quadd}{\quad \quad}

\newcommand{\wcset}{\what{\hspace{0pt}\mathcal{C}}}
\newcommand{\wpset}{\what{\hspace{0pt}\mathcal{P}}}
\newcommand{\sys}{\mathcal{C}}
\newcommand{\pset}{\mathcal{P}}
\newcommand{\what}[1]{\widehat{#1}}
\newcommand{\wtilde}[1]{\widetilde{#1}}
\newcommand{\wg}{\widehat{\gamma}}
\newcommand{\wc}{\widehat{c}}
\renewcommand{\wp}{\widehat{p}}
\newcommand{\wrho}{\widehat{\rho}}
\renewcommand{\th}[1]{\tilde{\hat{#1}}}

\newcommand{\aequation}[1]{\begin{equation} \begin{aligned} #1 \end{aligned} \end{equation}}

\newcommand{\calO}{\mathcal{O}}

\newtheorem{theorem}{Theorem} 
\newtheorem{proposition}[theorem]{Proposition}
\newtheorem{corollary}[theorem]{Corollary}
\newtheorem{lemma}[theorem]{Lemma}
\newtheorem{remark}[theorem]{Remark}

\begin{document}

\begin{frontmatter}

\title{An online-adaptive hyperreduced reduced basis element method for parameterized component-based nonlinear systems using hierarchical error estimation}

\author[affl1,affl2]{Mehran Ebrahimi\corref{cor1}}
\ead{m.ebrahimi@mail.utoronto.ca}
\author[affl1]{Masayuki Yano}
\ead{masa.yano@utoronto.ca}

\affiliation[affl1]{
  organization={Institute for Aerospace Studies, University of Toronto},
  addressline={4925 Dufferin Street},
  city={Toronto},
  postcode={M3H 5T6},
  state={Ontario},
  country={Canada}}
\affiliation[affl2]{
  organization={Autodesk Research},
  addressline={661 University Avenue},
  city={Toronto},
  postcode={M5G 1M1},
  state={Ontario},
  country={Canada}}
\cortext[cor1]{Corresponding author}

\begin{abstract}
\noindent
We present an online-adaptive hyperreduced reduced basis element method for model order reduction of parameterized, component-based nonlinear systems. The method, in the offline phase, prepares a library of hyperreduced archetype components of various fidelity levels and, in the online phase, assembles the target system using instantiated components whose fidelity is adaptively selected to satisfy a user-prescribed system-level error tolerance. To achieve this, we introduce a hierarchical error estimation framework that compares solutions at successive fidelity levels and drives a local refinement strategy based on component-wise error indicators. We also provide an efficient estimator for the system-level error to ensure that the adaptive strategy meets the desired accuracy. Component-wise hyperreduction is performed using an empirical quadrature procedure, with the training accuracy guided by the Brezzi--Rappaz--Raviart theorem. The proposed method is demonstrated on a family of nonlinear thermal fin systems comprising up to 225~components and 68 parameters. Numerical results show that the hyperreduced basis element model achieves $\calO(100)$ computational reduction at 1\% error level relative to the truth finite-element model. In addition, the adaptive refinement strategy provides more effective error control than uniform refinement by selectively enriching components with higher local errors.
\end{abstract}


\end{frontmatter}

\section{Introduction}
\label{sec:introduction}
Many scenarios in computational science, including design optimization, uncertainty quantification, and control, frequently involve many-query problems that require repeated solutions of parameterized partial differential equations (PDEs). For problems whose solution manifold admits accurate approximation in a low-dimensional space, reduced basis (RB) methods offer an effective approach to rapidly and reliably approximate the PDE solution at different parameter values~\cite{rozza2008reduced, quarteroni2015reduced, benner2015survey, hesthaven2016certified}. A typical workflow for RB methods involves separating the computation into offline and online phases. In the offline phase, the high-fidelity (i.e., \emph{truth}) problem is solved (using, for example, finite element (FE) methods) for many training parameter values to generate solution snapshots, which are then used to construct a basis for the RB space. For nonlinear problems, this phase also includes hyperreduction~\cite{barrault2004empirical, grepl2007efficient, hesthaven2016certified}. Although the offline phase can be computationally demanding, it enables significant computational savings in the online phase, where the reduced problem is solved numerous times for different parameter values in the many-query application.

Despite the effectiveness of standard RB methods for certain problems, they face significant challenges when applied to some complex engineering applications. One major limitation is their inability to accommodate topology-varying parameterizations: even slight changes in the domain topology can render a trained RB model inapplicable. Additionally, these methods are typically constrained to problems with a small number of parameters, as the offline phase requires repeated high-fidelity simulations, which become prohibitively expensive for large-scale problems. Consequently, the RB space constructed from a limited set of snapshots often lacks sufficient expressiveness to generalize to unseen parameters, which ultimately limits the robustness and scalability of the approach.

To overcome these challenges, \emph{component-based} RB methods have been developed~\cite{bourquin1992component, maday2002reduced}. These methods adopt a divide-and-conquer strategy by decomposing the global domain into smaller, computationally manageable subdomains. During the offline phase, a library of interoperable \emph{archetype} components and their associated local RB spaces is constructed. In the online phase, copies of the archetype components in the library are \emph{instantiated} to match the specific topological configuration of the system, and a global RB model for the entire system is assembled by coupling the preconstructed local reduced models. This approach eliminates the need to retrain for each new system configuration and facilitates reduced-order modeling of large-scale problems by avoiding the costly generation of global solution snapshots.

Component-based RB methods have been applied to both linear and nonlinear problems, although the majority of studies focus on linear problems. The reduced basis element (RBE) method~\cite{maday2002reduced,maday2004reduced,lovgren2006reduced} combines domain decomposition with RB methods and uses Lagrange multipliers to couple local reduced models in the online phase. The static condensation RBE method~\cite{huynh2013static, huynh2013staticc, smetana2016optimal} decomposes the degrees of freedom (DoFs) in each component into \emph{port} (interface) and \emph{bubble} (interior) DoFs. It employs static condensation~\cite{wilson1974static} to form a Schur complement system involving only port DoFs and applies RB approximations within each component to reduce computational cost and accommodate parametric variations. A port-reduced variant of this method~\cite{eftang2013port, eftang2014port, smetana2015new} further reduces the size of the Schur complement system by approximating the solution on global ports through RB methods applied to port modes. Hoang et al.~\cite{hoang2021domain} introduce a non-overlapping domain-decomposition least-squares Petrov--Galerkin method that weakly enforces interface continuity between subdomains via compatibility constraints. Iollo et al.~\cite{iollo2023one} develop a component-based model reduction approach for parameterized nonlinear elliptic PDEs that employs overlapping subdomains and an optimization-based formulation to minimize solution jumps across component interfaces. Smetana and Taddei~\cite{smetana2023localized} propose an overlapping multidomain RB method that uses the partition-of-unity method. Diaz et al.~\cite{diaz2024fast} integrate nonlinear approximation spaces, generated via autoencoders, with domain decomposition to address problems with slowly decaying Kolmogorov $n$-widths. In~\cite{ebrahimi2024hyperreduced}, we develop a hyperreduced RBE (HRBE) method that applies an online adaptive scheme, informed by the Brezzi--Rappaz--Raviart (BRR) theorem~\cite{brezzi1980finite, caloz1997numerical}, to select the appropriate hyperreduction fidelity for each component that ensures the satisfaction of the user-prescribed system-level error in the online phase. Finally, Chung et al.~\cite{chung2024train} devise a non-overlapping component-based RB method for linear problems using a discontinuous Galerkin domain decomposition and a physics-constrained, data-driven strategy, which is later extended in~\cite{chung2024scaled} to steady Navier--Stokes equations.

An important consideration in component-based RB methods is determining the appropriate dimension (i.e., fidelity) of local RB models. Overly large local RB dimensions increase the computational cost and memory footprint of the online phase, while excessively small dimensions compromise solution quality. The challenge is further compounded by the fact that, in component-based RB methods, the specific systems into which the trained archetype components will be integrated are unknown during the training phase. In this study, we extend our previous work~\cite{ebrahimi2024hyperreduced} by introducing an online-adaptive HRBE method to address this challenge for nonlinear problems. Additionally, we employ port reduction to further reduce the computational cost and memory requirements of the online phase. To ensure that the user-prescribed system-level error tolerance is met, we propose a refinement strategy that adaptively refines the RB fidelity for individual components and ports during the online phase. This strategy is based on a \emph{hierarchical} error estimation framework~\cite{hain2019hierarchical}, wherein two approximate solutions with different RB space fidelities are computed and compared at each refinement step. Specifically, a coarser solution is first computed, and its error is estimated by comparing it against a refined solution obtained with an enriched RB space. The refinement process proceeds iteratively until the prescribed accuracy threshold is satisfied. In the component-based setting, we \emph{locally} increase RB fidelities only where necessary to maintain computational efficiency.

Our adaptive refinement strategy shares similarities with the algorithms introduced in~\cite{ohlberger2015error, buhr2018exponential, smetana2023localized} but differs in several key aspects. In~\cite{ohlberger2015error, buhr2018exponential}, the accuracy of local approximation spaces is improved during the online phase by enriching them with local truth snapshot solutions. In~\cite{smetana2023localized}, local RB spaces are enriched during the offline phase using global reduced solves and a local residual-based error indicator; however, the algorithm in~\cite{smetana2023localized} relies on computing the dual norm of local truth residuals and does not incorporate hyperreduction. Our approach, in contrast, does not depend on any truth quantities during the online phase and further addresses the implications of adaptivity in component-wise hyperreduction training during the offline phase.

The contributions of the present work are sixfold:
\begin{enumerate}
  \item We develop an online-adaptive HRBE method that enables model order reduction of parameterized, component-based nonlinear systems. In the offline phase, the method constructs a library of multi-fidelity hyperreduced components. In the online phase, the method adaptively selects the appropriate fidelity of components and ports in the system to ensure that the user-prescribed system-level error tolerance is met efficiently.
  \item We devise a port reduction strategy to obtain compact modal representations of solutions at component interfaces, which leads to additional computational savings and reduced memory requirements during the online phase.
  \item We appeal to the BRR theorem to develop a hyperreduction-fidelity selection mechanism for the component-wise hyperreduction procedure introduced in~\cite{ebrahimi2024hyperreduced}.
  \item We introduce an online-efficient, hierarchical error estimation framework to compute the system-level error estimate and component-wise error indicators during the online phase.
  \item We develop an adaptive refinement algorithm informed by the component-wise error indicators to selectively enrich the RB and hyperreduction fidelities of components and ports during the online phase.
  \item We demonstrate the effectiveness of the proposed online-adaptive HRBE method on a family of nonlinear thermal fin systems with up to 225 instantiated components and 68 independent parameters. 
\end{enumerate}

The remainder of the paper is organized as follows. Section~\ref{sec:def} introduces the model problem and defines the notions of components, ports, and system. Section~\ref{sec:hrbe} presents the port-reduced HRBE method, where we apply the bubble--port decomposition of functions and derive the truth, RB, and HRBE problems in terms of bubble and port solutions. Section~\ref{sec:adaptive} describes the adaptive HRBE method, which uses a library of multi-fidelity hyperreduced components and a hierarchical error estimator. Section~\ref{sec:training} details the component-wise offline training procedures. Section~\ref{sec:error} establishes the theoretical foundation for fidelity selection in the component-wise hyperreduction. Section~\ref{sec:res} presents numerical experiments that demonstrate the effectiveness of the proposed online-adaptive HRBE method. Finally, Section~\ref{sec:conclusion} concludes the paper.

\section{Model problem}
\label{sec:def}
In this section, we present the general form of the considered model problem. To maintain consistency, we adopt the notation introduced in our earlier work~\cite{ebrahimi2024hyperreduced}, with slight modifications to incorporate port-reduction.

\subsection{Components, ports, and system}
\label{subsec:comps}
We begin by introducing the entities associated with archetype components and ports. We define a library of $\what{N}_\comp$ \emph{archetype components} and $\what{N}_\port$ \emph{archetype ports}\footnote{Throughout this document, the notation $\what{\cdot}$ refers to quantities associated with or defined over the archetype (as opposed to instantiated) components and ports.}. We introduce $\wcset$ and $\wpset$ as the set of archetype components and ports in the library, respectively. For each archetype component $\what{c} \in \wcset$, we introduce $\what{\Omega}_{\what{c}} \subset \mathbb{R}^d$, $\partial \what{\Omega}_{\what{c}} \subset \mathbb{R}^{d}$, ${\what{\mathcal{D}}}_{\what{c}} \subset \mathbb{R}^{{n}_{\what{c}}}$, and ${\what{\mu}}_{\what{c}} \in {\what{\mathcal{D}}}_{\what{c} }$ as, respectively, its bounded $d$-dimensional reference spatial domain, its Lipschitz-continuous boundary, its bounded ${n}_{\what{c}}$-dimensional parameter domain, and ${n}_{\what{c}}$-tuple specifying its reference parameter values. The parameter domain ${\what{\mathcal{D}}}_{\what{c}}$ may include both \emph{geometric} and \emph{non-geometric} parameters. Geometric parameters describe variations in the component's shape, while non-geometric parameters account for attributes that do not affect the shape, such as material properties or traction in elasticity applications.

\begin{figure}[!bt]
\centering
  \begin{subfigure}[b]{0.3\textwidth}
    \centering
      \includegraphics[width=0.23\textwidth]{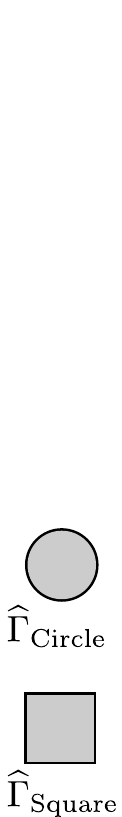}
      \caption{Archetype ports.}\label{fig:sysexamplea}
  \end{subfigure}
  \hspace{-1.5em}
  \begin{subfigure}[b]{0.3\textwidth}
      \includegraphics[width=\textwidth]{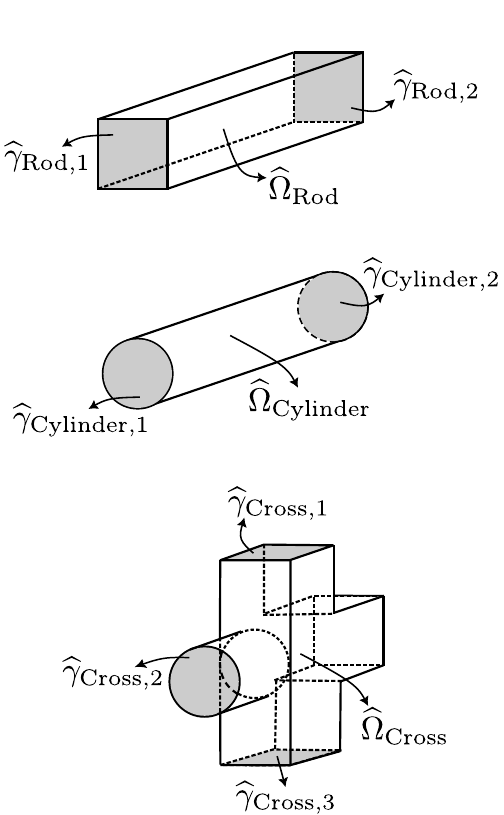}
      \caption{Archetype components.}\label{fig:sysexampleb}
  \end{subfigure}
  \hspace{0em}
    \begin{subfigure}[b]{0.3\textwidth}
        \includegraphics[width=\textwidth]{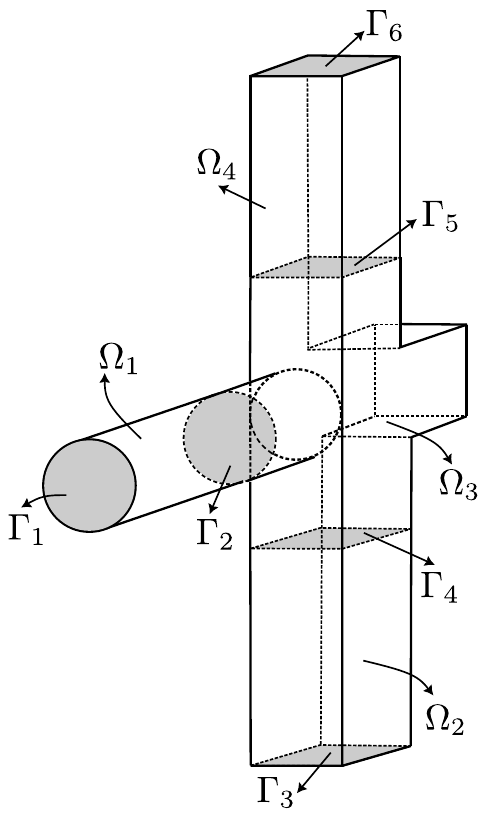}
        \caption{Assembled system.}\label{fig:sysexamplec}
    \end{subfigure}
    \caption{(a) Archetype ports, (b) archetype components with ports mapped from the archetype ports in (a), and (c) a system with four instantiated components and six global ports.}
    \label{fig:sysexample}
\end{figure}

Each archetype component $\what{c} \in \wcset$ has $n_{\what{c}}^{\gamma}$ disjoint \emph{local ports}, created from the reference archetype ports in $\wpset$ through a geometric mapping. We introduce $\what{\Gamma}_{\wp} \subset \mathbb{R}^{d-1}$ and $\what{\gamma}_{\what{c},p} \subset \mathbb{R}^{d}$ as, respectively, the domain of $\wp \in \wpset$ and $p \in \pset_{\wc} \equiv \{ 1, \dots, n_{\wc}^{\gamma} \}$. We further introduce $\mathcal{R}_{\wc, p}: \what{\Gamma}_{\pi_\wc(p)} \to \what{\gamma}_{\what{c},p}$ as the invertible, geometric-parameter-\emph{independent} mapping between $\what{\gamma}_{\what{c},p}$ and $\what{\Gamma}_{\pi_\wc(p)}$ such that $\what{\gamma}_{\what{c},p} = \mathcal{R}_{\wc, p}(\what{\Gamma}_{\pi_\wc(p)})$, where $\pi_\wc: \pset_\wc \rightarrow \wpset$ is a map from the local ports of $\wc$ to their corresponding (exactly one) archetype port in $\wpset$. We assume the boundary of all components is Lipschitz-continuous and all ports of an archetype component are mutually separated by a non-port boundary surface. Figure~\ref{fig:sysexampleb} shows three archetype components whose ports are mapped from the archetype ports in Figure~\ref{fig:sysexamplea}.

We next introduce the entities associated with instantiated components. We define $\sys$ as a set of $N_\comp$ instantiated components that create a system. Each instantiated component is generated from an archetype component in the library through a (parameterized) geometric mapping. We introduce $M: \sys \rightarrow \wcset$ as a map from the instantiated components to their corresponding (exactly one) archetype component in the library. We further define $\Omega_c \subset \mathbb{R}^d$ and $\mu_c \in \mathcal{D}_c \equiv \widehat{\mathcal{D}}_{{M}(c)}$ as the physical domain of the instantiated component $c \in \sys$ and its parameter tuple, respectively. The parameterized geometric mappings between archetype and instantiated component domains are denoted by $\mathcal{G}_c: \widehat{\Omega}_{{M}(c)} \times \mathcal{D}_{c} \rightarrow \Omega_c$ such that $\Omega_c = \mathcal{G}_c (\widehat{\Omega}_{{M}(c)}; \mu_c)$. The physical domain of the $p$-th local port of $c$, with $\: p \in {\pset}_{M(c)}$, is given by $\gamma_{c,p} \equiv \mathcal{G}_c (\widehat{\gamma}_{{M}(c),p}; \mu_c) = \mathcal{G}_c (\mathcal{R}_{M(c),p}(\what{\Gamma}_{\pi_{M(c)}(p)}); \mu_c)$.

The components in the system are connected through their local ports, which results in $N_\port$ \emph{global ports}. We denote the set of global ports by $\pset \equiv \{1,\dots, N_\port\}$. Geometric mappings ensure that the ports conform to one another. Each global port is assumed to be shared by at most two instantiated components. A local port on the system boundary also forms a global port. Essential boundary conditions at the system level are applied to these boundary global ports. Figure~\ref{fig:sysexamplec} illustrates an example of a system composed of four components and six global ports. 

\subsection{Exact problem formulation}
\label{subsec:continuous}
We begin by introducing the function spaces associated with archetype and instantiated components. For each archetype component ${\what{c}} \in \wcset$, we define a Hilbert space $\widehat{\mathcal{V}}_{{\what{c}}} \subset H^1(\widehat{\Omega}_{\what{c}})$, equipped with an inner product $(\cdot, \cdot)_{\widehat{\mathcal{V}}_{{\what{c}}}}$ and the induced norm $|| \cdot ||_{\widehat{\mathcal{V}}_{{\what{c}}}} \equiv \sqrt{(\cdot, \cdot)_{\widehat{\mathcal{V}}_{{\what{c}}}}}$, equivalent to the $H^1(\widehat{\Omega}_{\what{c}})$ norm. For each instantiated component $c \in \sys$, we introduce the geometric-parameter-\emph{dependent} space $\mathcal{V}_{c} \equiv \left\{v = \widehat{v} \circ \mathcal{G}_c^{-1}(\cdot; \mu_c) : \widehat{v} \in \widehat{\mathcal{V}}_{{M}(c)} \right\}$, along with its inner product and induced norm $(\cdot, \cdot)_{\mathcal{V}_c}$ and $\|\cdot \|_{\mathcal{V}_c} \equiv \sqrt{(\cdot, \cdot)_{\mathcal{V}_c}}$, respectively.

We present the \emph{exact} nonlinear model problem in its weak form as follows: given $\mu \equiv (\mu_c)_{c \in \sys} \in \mathcal{D} \equiv \prod_{c \in \sys} \mathcal{D}_{c}$, find $u(\mu) \in\ \mathcal{V}$ such that
\begin{equation}\label{eq:exact}
  \begin{aligned}
  R(u(\mu), v;\mu) \equiv \sum_{c \in \sys}  {R}_{c} \left( u(\mu) \big|_{\Omega_c}, v \big|_{\Omega_c}; \mu_c  \right) = 0 \quad \forall v \in \mathcal{V},
  \end{aligned}
\end{equation}
where $\mathcal{V} = \left\{ v \in H^1 (\Omega) : v |_{\Gamma_D} = 0 \right\}$ and is endowed with the $H^1({\Omega})$ norm $\|\cdot \|_{\mathcal{V}} \equiv \sqrt{\sum_{c \in \sys} \|\cdot \|_{\mathcal{V}_c}^2}$, ${R}_c : {\mathcal{V}}_c \times {\mathcal{V}}_c \times {\mathcal{D}}_c \rightarrow \mathbb{R}$ is the physical-domain residual form of the instantiated component $c \in \sys$, $\Omega$ is the system's physical domain such that $\overline \Omega = \cup_{c \in \sys} \overline \Omega_c$, and $\Gamma_D$ is its Dirichlet boundary, which is assumed non-empty.\footnote{To streamline the presentation, we assume homogeneous Dirichlet boundary conditions; problems involving nonhomogeneous boundary conditions can be handled with minor modifications.} The residual form $R_c(\cdot,\cdot;\cdot)$ takes on the form 
\begin{equation}\label{eq:integrand}
  \begin{aligned}
  {R}_c (w,v;\mu) &= \int_{{\Omega}_c} {r}_c(w,v;{x},\mu)  \: d{x} & &\quad \forall w,v \in {\mathcal{V}}_c, \: \forall \mu \in {\mathcal{D}}_c,
  \end{aligned}
\end{equation}
where ${r}_c: {\mathcal{V}}_c \times {\mathcal{V}}_c \times {\Omega}_c \times {\mathcal{D}}_c \rightarrow \mathbb{R}$ is the physical-domain integrand, which is linear in its second argument but, in general, nonlinear in its first argument; for example, for the nonlinear heat equation considered in Section~\ref{sec:res}, the integrand is given by $r_c(w,v;{x},\mu) = \nabla v \cdot k(w) \nabla w - v f(\mu)$ for some nonlinear diffusion coefficient $k(w)$ and heat source $f(\mu)$. We assume that the exact problem~\eqref{eq:exact} is well-posed $\forall \mu \in \mathcal{D}$. 

To facilitate the treatment of geometric parameters, we express the system-level residual in terms of the reference domain of the components as
\begin{equation}\nonumber
  \begin{aligned}
  R(w, v;\mu) &= \sum_{c \in \sys}  \widehat{R}_{M(c)} \left( w \big|_{\Omega_c} \circ \mathcal{G}_c(\cdot; \mu_c), v \big|_{\Omega_c} \circ \mathcal{G}_c(\cdot; \mu_c); \mu_c  \right) \quadd \forall w,v \in \mathcal{V}, \: \forall \mu \in \mathcal{D}.
  \end{aligned}
\end{equation}
Here, $\widehat{R}_{\what{c}} : \widehat{\mathcal{V}}_{\what{c}} \times \widehat{\mathcal{V}}_{\what{c}} \times \widehat{\mathcal{D}}_{\what{c}} \rightarrow \mathbb{R}$ is the reference-domain residual of the archetype component $\what{c} \in \wcset$ given by
\begin{equation}\nonumber
\begin{aligned}
\widehat{R}_\wc (w,v;\mu) &= \int_{\widehat{\Omega}_\wc} \widehat{r}_\wc (w,v;\what{x},\mu)  \: d\what{x} & &\quad \forall w,v \in \widehat{\mathcal{V}}_\wc, \: \forall \mu \in \widehat{\mathcal{D}}_\wc,
\end{aligned}
\end{equation}
where $\widehat{r}_\wc: \widehat{\mathcal{V}}_\wc \times \widehat{\mathcal{V}}_\wc \times \widehat{\Omega}_\wc \times \widehat{\mathcal{D}}_\wc \rightarrow \mathbb{R}$ is the reference-domain integrand. For each $c \in \sys$, the physical-domain integrand $r_c(\cdot,\cdot;\cdot,\cdot)$ satisfies
\begin{equation}\nonumber
    \begin{aligned}
    {r}_c(w,v;{x},\mu_c) &= \what{r}_{M(c)}(w \circ \mathcal{G}_c(\cdot; \mu_c),v \circ \mathcal{G}_c(\cdot; \mu_c); \mathcal{G}_c^{-1}(x; \mu_c),\mu_c)  \: \det(\mathcal{J}_c(\mathcal{G}_c^{-1}(x; \mu_c); \mu_c))^{-1}
\end{aligned}
\end{equation}
for all $w,v \in \mathcal{V}_c$, $x \in \Omega_c$, and $\mu_c \in \mathcal{D}_c$, where $\mathcal{J}_c(\cdot; \mu_c): \widehat{\Omega}_{\wc} \to \mathbb{R}^{d \times d}$ is the Jacobian of $\mathcal{G}_c(\cdot; \mu_c)$.
\section{Port-reduced HRBE method}
\label{sec:hrbe}
In this section, we present a port-reduced HRBE method to accurately approximate the truth solution while significantly reducing computational cost. Computational efficiency is achieved through the application of (i) RB approximations to the bubble and port spaces of the components and (ii) component-wise hyperreduction. To this end, we first introduce the bubble--port decomposition of functions and formulate the truth problem in terms of bubble and port functions. We then present the port-reduced HRBE method.

\subsection{Bubble--port decomposition of functions}
We follow the decomposition approach introduced in~\cite{huynh2013static}. We begin by defining approximation spaces for archetype ports. For each $\wp \in \wpset$, we introduce a Hilbert space $\widehat{\mathcal{X}}_{\wp} \subset H^1(\widehat{\Gamma}_{\wp})$ endowed with an inner product $(\cdot, \cdot)_{\widehat{\mathcal{X}}_{\wp}}$ and the induced norm $|| \cdot ||_{\widehat{\mathcal{X}}_{{\wp}}} \equiv \sqrt{(\cdot, \cdot)_{\widehat{\mathcal{X}}_{{\wp}}}}$, equivalent to the $H^1(\widehat{\Gamma}_{\wp})$ norm. For this port, we introduce an $\mathcal{N}_{\wp}$-dimensional FE space $\widehat{\mathcal{X}}_{h,{\wp}} \subset \widehat{\mathcal{X}}_{\wp}$ constructed by tessellating $\what{\Gamma}_\wp$ with nonoverlapping, conforming elements. For each $\wp$, we introduce $\{ \widehat{\tau}_{\wp,i} \}_{i=1}^{\mathcal{N}_{\wp}}$ as the eigenbasis associated with an eigenproblem: find eigenpairs $(\widehat{\tau}_{\wp,i}, \lambda_{\wp,i}) \in \widehat{\mathcal{X}}_{h,\wp} \times \mathbb{R}$, $i = 1,\dots,\mathcal{N}_{\wp}$, such that
\begin{equation}\nonumber
  \begin{aligned}
  \int_{\widehat{\Gamma}_{\wp}} \nabla \widehat{\tau}_{\wp,i} \cdot \nabla \what{y} \: ds &= \lambda_{\wp,i} \int_{\widehat{\Gamma}_{\wp}} \widehat{\tau}_{\wp,i} \: \what{y} \: ds \quad \forall \what{y} \in \widehat{\mathcal{X}}_{h,\wp}, \\
  \left\| \widehat{\tau}_{\wp,i} \right\|_{L^2(\widehat{\Gamma}_{\wp})} &= 1.
  \end{aligned}
\end{equation}

Next, we define approximation spaces associated with archetype components. For each $\wc \in \wcset$, we introduce an $\mathcal{N}_{\what{c}}^\bb$-dimensional bubble FE space $\widehat{\mathcal{V}}_{h,{\what{c}}}^\bb \equiv \left\{ v \in \widehat{\mathcal{V}}_{h,{\what{c}}} : \: v \big|_{\wg_{\wc,p}} = 0 \ \forall p \in {\pset}_{\what{c}} \right\}$, where $\widehat{\mathcal{V}}_{h,{\what{c}}} \subset \widehat{\mathcal{V}}_{{\what{c}}}$ is an $\mathcal{N}_{\what{c}}$-dimensional FE space formed via a tessellation of $\widehat{\Omega}_{\what{c}}$ into nonoverlapping, conforming elements. Additionally, for each $p \in \pset_\wc$, we introduce an $\mathcal{N}_{\pi_\wc(p)}$-dimensional port FE space $\widehat{\mathcal{X}}_{h,{\what{c}}}^{p} \equiv \left\{v = \widehat{v} \circ \mathcal{R}_{\wc,p}^{-1}(\cdot) : \: \widehat{v} \in \widehat{\mathcal{X}}_{h,\pi_\wc(p)} \right\}$ with the basis $\{ \widehat{\tau}_{\wc,i}^p \equiv \what{\tau}_{\pi_{\wc}(p),i} \circ \mathcal{R}_{\wc,p}^{-1}(\cdot) \}_{i=1}^{\mathcal{N}_{\pi_\wc(p)}}$. We elliptically lift these basis functions to the interior of $\wc$ to obtain $\{ \widehat{\psi}_{\wc,i}^p \in \widehat{\mathcal{V}}_{h,\wc} \}_{i=1}^{\mathcal{N}_{\pi_\wc(p)}}$ such that
\begin{equation}\nonumber
  \begin{alignedat}{3}
  \int_{\widehat{\Omega}_\wc} \nabla \widehat{\psi}_{\wc,i}^p \cdot \nabla v \:  d\widehat{x} &= 0  \quad \quad && \forall v \in \widehat{\mathcal{V}}_{h,\wc}^\bb, \\
  \widehat{\psi}_{\wc,i}^p &= \widehat{\tau}_{\wc,i}^p \quad \quad && \text{on } \widehat{\gamma}_{\wc,p}, \\
  \widehat{\psi}_{\wc,i}^p &= 0 \quad \quad && \text{on } \widehat{\gamma}_{\wc,p'} \quad \forall p' \neq p,
  \end{alignedat}
\end{equation}
for all $p \in \pset_\wc$ and define $\what{\mathcal{V}}^p_{h,\wc} \equiv \text{span} \{ \what{\psi}^p_{\wc, i} \}_{i=1}^{\mathcal{N}_{\pi_{\wc}(p)}}$. We note that the approximation spaces satisfy $\widehat{\mathcal{V}}_{h,{\what{c}}}|_{\wg_{\wc,p}} = \what{\mathcal{V}}^p_{h,\wc}|_{\wg_{\wc,p}}$ $\forall \wc \in \wcset$ and $\forall p \in \pset_\wc$.

We can now present the bubble--port decomposition of functions defined on each instantiated component $c \in \sys$. We introduce FE spaces for (full) component, component bubble, port modes, lifted port modes, and their collection:
\begin{equation}
    \begin{aligned}
        \mathcal{V}_{h,c} &\equiv \left\{v = \widehat{v} \circ \mathcal{G}_c^{-1}(\cdot; \mu_c) : \: \widehat{v} \in \widehat{\mathcal{V}}_{h, {M}(c)} \right\} \subset \mathcal{V}_{c}, \nonumber \\
        \mathcal{V}_{h,c}^\bb &\equiv \Big\{v = \widehat{v} \circ \mathcal{G}_c^{-1}(\cdot; \mu_c) :\: \widehat{v} \in \widehat{\mathcal{V}}_{h, {M}(c)}^\bb \Big\} \subset \mathcal{V}_{h,c}, \nonumber \\
        \mathcal{X}_{h,c}^p &\equiv \Big\{v = \widehat{v} \circ \mathcal{G}_c^{-1}(\cdot; \mu_c) :\: \widehat{v} \in \widehat{\mathcal{X}}_{h, {M}(c)}^p \Big\} \quad \forall p \in \pset_{M(c)}, \nonumber \\
        \mathcal{V}_{h,c}^p &\equiv \left\{v = \widehat{v} \circ \mathcal{G}_c^{-1}(\cdot; \mu_c) : \: \widehat{v} \in \widehat{\mathcal{V}}_{h,M(c)}^p \right\} \quad \forall p \in \pset_{M(c)}, \nonumber \\
        \mathcal{V}_{h,c}^\gamma &\equiv \cup_{p \in \pset_{M(c)}} \mathcal{V}_{h,c}^p.
    \end{aligned}
\end{equation}
Therefore, any $v_{h,c} \in \mathcal{V}_{h,c}$ can be written as 
\begin{equation}\label{eq:compsol}
    v_{h,c} = v^\bb_{h,c} + v^\gamma_{h,c} = v^\bb_{h,c} + \sum_{p\in \pset_{M(c)}} v^p_{h,c},
\end{equation}
where $v^\bb_{h,c} \in \mathcal{V}_{h,c}^\bb$ and $v^\gamma_{h,c} \in \mathcal{V}_{h,c}^\gamma$ are, respectively, the bubble and port parts of $v_{h,c}$.

\subsection{Truth problem formulation}
We now formulate the truth problem in terms of bubble and port functions. We define $(\widehat{x}_{\wc,q}, \wrho_{\wc,q} )_{q=1}^{Q_\wc}$ $\forall \wc \in \wcset$ as the truth quadrature rule in the reference domain $\widehat{\Omega}_\wc$ of each archetype component $\wc \in \wcset$. The truth problem is the following: given $\mu = (\mu_c)_{c \in \sys} \in \mathcal{D}$, find $\{ {u}_{h,c}^{\bb}(\mu) \in \mathcal{V}_{h,c}^\bb \}_{c \in \sys}$ and $\{ {u}_{h,c}^\gamma(\mu) \in \mathcal{V}_{h,c}^\gamma \}_{c \in \sys}$ such that, for all $\{ v_{h,c}^{\bb} \in \mathcal{V}_{h,c}^\bb \}_{c \in \sys}$ and $\{ {v}_{h,c}^\gamma \in \mathcal{V}_{h,c}^\gamma \}_{c \in \sys}$,
\begin{equation}\label{eq:truthref}
  \begin{aligned}
  R_h(u_h(\mu), v_h;\mu) &\equiv \sum_{c \in \sys} \sum_{q=1}^{Q_{M(c)}} \wrho_{M(c),q} \widehat{r}_{M(c)} \biggl( \Bigl[ u^\bb_{h,c}(\mu) + u^{\gamma}_{h,c}(\mu) \Bigr] \circ \mathcal{G}_{c}(\widehat{x}_{M(c),q}; \mu_c), \\
  &\hspace{11.5em}\Bigl[ v^\bb_{h,c} + v^{\gamma}_{h,c} \Bigr] \circ \mathcal{G}_{c}(\widehat{x}_{M(c),q}; \mu_c); \widehat{x}_{M(c),q},\mu_c \biggr) = 0,
\end{aligned}
\end{equation}
so that the system-level truth solution $u_h(\mu) \in \mathcal{V}_h$ is given by 
\begin{equation}\nonumber
  u_h(\mu) = \sum_{c \in \sys} \Big[ u^\bb_{h,c}(\mu) + u^{\gamma}_{h,c}(\mu) \Big] = \sum_{c \in \sys} \Big[ u^\bb_{h,c}(\mu) + \sum_{p \in \pset_{M(c)}} u^{p}_{h,c}(\mu) \Big].
\end{equation}
Here, $\mathcal{V}_{h} = \left( \boplus_{c \in \sys}\mathcal{V}_{h,c} \right) \cap \mathcal{V}$ is the $\mathcal{N}_h$-dimensional system-level truth FE space, where the intersection with $\mathcal{V}$ enforces the Dirichlet boundary conditions and continuity at the global ports. Similar to the exact problem in~\eqref{eq:exact}, we assume that the truth problem is well-posed $\forall \mu \in \mathcal{D}$. As mentioned earlier, we assume the connected local ports at the system level are conformal. Thus, for the $p$-th global port, $p \in \pset$, shared by the $l$-th port of $c \in \sys$ and the $l'$-th port of $c' \in \sys$, we have $\mathcal{X}_{h,c}^l = \mathcal{X}_{h,c'}^{l'}$ and $u_{h,c}(\mu)\big|_{\gamma_{c,l}} = u_{h,c'}(\mu)\big|_{\gamma_{c',l'}}$.

\subsection{Port-reduced RB and HRBE problem formulations}
\label{subsec:hrbe}
We begin by introducing the bubble and port RB spaces. We assume that $\forall \wc \in \wcset$ and $\forall \wp \in \wpset$, the parametric manifold spanned by their respective truth solutions due to solving~\eqref{eq:truthref} for all $\mu \in \mathcal{D}$ are amenable to accurate approximation by a low-dimensional linear space. For each $\wp \in \wpset$, we introduce an $N_\wp \ll \mathcal{N}_\wp$-dimensional RB space $\what{\mathcal{X}}_{\rb,\wp} \subset \what{\mathcal{X}}_{h,\wp}$. We further introduce, for each $\wc \in \wcset$, an $N^\bb_\wc \ll \mathcal{N}^\bb_\wc$-dimensional RB space $\what{\mathcal{V}}^\bb_{\rb,\wc} \subset \what{\mathcal{V}}^\bb_{h,\wc}$. We denote the basis of $\what{\mathcal{X}}_{\rb,\wp}$ and $\what{\mathcal{V}}^\bb_{\rb,\wc}$ by $\{\what{\chi}_{\wp,i} \}_{i=1}^{N_{\wp}}$ and $\{\what{\xi}^\bb_{\wc,i} \}_{i=1}^{N^\bb_{\wc}}$, respectively. The computational procedures for constructing these RB spaces are discussed in Section~\ref{sec:training}; for now, we assume $\{\what{\chi}_{\wp,i} \}_{i=1}^{N_{\wp}}$ and $\{\what{\xi}^\bb_{\wc,i} \}_{i=1}^{N^\bb_{\wc}}$ are given. Subsequently, for all $p \in \pset_\wc$, we introduce the $N_{\pi_\wc(p)}$-dimensional port RB space $\widehat{\mathcal{X}}_{\rb,{\what{c}}}^{p} \equiv \left\{v = \widehat{v} \circ \mathcal{R}_{\wc,p}^{-1}(\cdot) : \: \widehat{v} \in \widehat{\mathcal{X}}_{\rb,\pi_\wc(p)} \right\} \subset \widehat{\mathcal{X}}_{h,{\what{c}}}^{p}$ spanned by $\{ \what{\chi}^p_{\wc,i} \equiv \what{\chi}_{\pi_{\wc}(p),i} \circ \mathcal{R}_{\wc,p}^{-1}(\cdot) \}_{i=1}^{N_{\pi_\wc(p)}}$. We also introduce $\{ \what{\theta}^p_{\wc,i} \in \widehat{\mathcal{V}}_{h,\wc}\}_{i=1}^{N_{\pi_{\wc}(p)}}$ $\forall p \in \pset_\wc$ obtained through elliptically lifting $\{ \what{\chi}^p_{\wc,i} \}_{i=1}^{N_{\pi_{\wc}(p)}}$ such that
\begin{equation}\label{eq:portlift}
  \begin{alignedat}{3}
  \int_{\widehat{\Omega}_\wc} \nabla \widehat{\theta}_{\wc,i}^p \cdot \nabla v \:  d\widehat{x} &= 0  \quad \quad && \forall v \in \widehat{\mathcal{V}}_{\rb,\wc}^\bb, \\
  \widehat{\theta}_{\wc,i}^p &= \widehat{\chi}_{\wc,i}^p \quad \quad && \text{on } \widehat{\gamma}_{\wc,p}, \\
  \widehat{\theta}_{\wc,i}^p &= 0 \quad \quad && \text{on } \widehat{\gamma}_{\wc,p'} \quad \forall p' \neq p,
\end{alignedat}
\end{equation}
and define $\what{\mathcal{V}}^p_{\rb,\wc} \equiv \text{span} \{ \what{\theta}^p_{\wc, i} \}_{i=1}^{N_{\pi_{\wc}(p)}}$. Finally, for any $\wc \in \wcset$, we introduce the $N_\wc$-dimensional RB space $\what{\mathcal{V}}_{\rb,\wc} \equiv \what{\mathcal{V}}_{\rb,\wc}^\bb \cup \left(\cup_{p \in \pset_\wc} \what{\mathcal{V}}^p_{\rb,\wc} \right)$, where $N_\wc = N_\wc^\bb + \sum_{p \in \pset_\wc} N_{\pi_\wc(p)} \ll \mathcal{N}_\wc$.

We now define the RB spaces for the instantiated components. For each $c \in \sys$, we introduce
\begin{equation}\label{eq:rbspaces}
    \begin{aligned}
      \mathcal{V}_{\rb,c}^\bb &\equiv \left\{v = \widehat{v} \circ \mathcal{G}_c^{-1}(\cdot; \mu_c) : \: \widehat{v} \in \widehat{\mathcal{V}}_{\rb,M(c)}^\bb \right\} \subset \mathcal{V}_{h,c}^\bb, \\
      \mathcal{V}_{\rb,c}^p &\equiv \left\{v = \widehat{v} \circ \mathcal{G}_c^{-1}(\cdot; \mu_c) : \: \widehat{v} \in \widehat{\mathcal{V}}_{\rb,M(c)}^p \right\} \subset \mathcal{V}_{h,c}^p \quad \forall p \in \pset_{M(c)}, \\
      \mathcal{V}_{\rb,c}^\gamma &\equiv \cup_{p \in \pset_{M(c)}} \mathcal{V}^p_{\rb,c} \subset \mathcal{V}^\gamma_{h,c}, \\
      \mathcal{V}_{\rb,c} &\equiv \mathcal{V}_{\rb,c}^\bb \cup \mathcal{V}_{h,c}^\gamma \subset \mathcal{V}_{h,c}.
    \end{aligned}
\end{equation}
Hence, any $v_{\rb,c} \in\mathcal{V}_{\rb,c}$ can be expressed as $v_{\rb,c} = v_{\rb,c}^\bb + v_{\rb,c}^\gamma$, where $v_{\rb,c}^\bb \in \mathcal{V}_{\rb,c}^\bb$ and $v_{\rb,c}^\gamma \in \mathcal{V}_{\rb,c}^\gamma$.

We now formulate the port-reduced RB problem: given $\mu = (\mu_c)_{c \in \sys} \in \mathcal{D}$, find $\{ {u}_{\rb,c}^{\bb}(\mu) \in \mathcal{V}_{\rb,c}^\bb \}_{c \in \sys}$ and $\{ {u}_{\rb,c}^\gamma(\mu) \in \mathcal{V}_{\rb,c}^\gamma\}_{c \in \sys}$ such that, for all $\{ v_{\rb,c}^{\bb} \in \mathcal{V}_{\rb,c}^\bb \}_{c \in \sys}$ and $\{ {v}_{\rb,c}^\gamma \in \mathcal{V}_{\rb,c}^\gamma \}_{c \in \sys}$,
\begin{equation}\label{eq:nlrbpr}
\begin{aligned}
{R}_{h} ({u}_{\rb}(\mu), {v}_{\rb};\mu ) = \sum_{c \in \sys} \sum_{q=1}^{{Q}_{M(c)}} \what{\rho}_{M(c),q} \: \widehat{r}_{M(c)} \biggl( &\Bigl[{u}^{\bb}_{\rb,c}(\mu) + {u}^{\gamma}_{\rb,c}(\mu) \Bigr] \circ {\mathcal{G}}_{c}(\what{x}_{M(c),q}; \mu_c),  \\
&\Bigl[{v}^{\bb}_{\rb,c} + v^{\gamma}_{\rb,c} \Bigr] \circ {\mathcal{G}}_{c}(\what{x}_{M(c),q}; \mu_c); \what{x}_{M(c),q},\mu_c \biggr) = 0,
\end{aligned}
\end{equation}
so that the system-level RB solution ${u}_\rb(\mu) \in \mathcal{V}_\rb$ is given by ${u}_\rb(\mu) = \sum_{c \in \sys} \left[ {u}^\bb_{\rb,c}(\mu) + {u}^{\gamma}_{\rb,c}(\mu) \right]$; here, $\mathcal{V}_{\rb} = \left( \boplus_{c \in \sys} \mathcal{V}_{\rb,c} \right) \cap \mathcal{V}$ is the $N_\rb$-dimensional system-level RB space. We assume that the port-reduced RB problem is well-posed $\forall \mu \in \mathcal{D}$.

We now present the port-reduced HRBE problem formulation. An essential component in formulating the problem is hyperreduction. To this end, we employ the EQP~\cite{Patera_2017_EQP, yano2019lp}, and more specifically its component-wise variant~\cite{ebrahimi2024hyperreduced}. Through the component-wise EQP, a sparse subset of reduced quadrature (RQ) points with re-weighted quadrature weights is found such that the integrals in the component residual forms are approximated to a prescribed accuracy. For each archetype component $\wc \in \wcset$, we introduce the residual RQ rule $(\th{x}_{\wc,q}, \th{\rho}_{\wc,q} )_{q=1}^{\wtilde{Q}_\wc} \subset (\widehat{x}_{\wc,q}, \wrho_{\wc,q} )_{q=1}^{Q_\wc}$, where $\wtilde{Q}_\wc \ll Q_\wc$. (More details on the component-wise EQP are provided in Section~\ref{sec:training}; for now, we assume $(\th{x}_{\wc,q}, \th{\rho}_{\wc,q} )_{q=1}^{\wtilde{Q}_\wc}$ for any $\wc \in \wcset$ is given.) We now present the port-reduced HRBE problem: given $\mu = (\mu_c)_{c \in \sys} \in \mathcal{D}$, find $\{ \wtilde{u}_{\rb,c}^{\bb}(\mu) \in \mathcal{V}_{\rb,c}^\bb \}_{c \in \sys}$ and $\{ \wtilde{u}_{h,c}^\gamma(\mu) \in \mathcal{V}_{\rb,c}^\gamma\}_{c \in \sys}$ such that, for all $\{ v_{\rb,c}^{\bb} \in \mathcal{V}_{\rb,c}^\bb \}_{c \in \sys}$ and $\{ {v}_{\rb,c}^\gamma \in \mathcal{V}_{\rb,c}^\gamma \}_{c \in \sys}$
\begin{equation}\label{eq:nlhrbepr}
\begin{aligned}
\wtilde{R}_{\rb} (\wtilde{u}_{\rb}(\mu), {v}_{\rb};\mu ) \equiv \sum_{c \in \sys} \sum_{q=1}^{\wtilde{Q}_{M(c)}} \th{\rho}_{M(c),q} \: \widehat{r}_{M(c)} \biggl( &\Bigl[ \wtilde{u}^{\bb}_{\rb,c}(\mu) + \wtilde{u}^{\gamma}_{\rb,c}(\mu) \Bigr] \circ {\mathcal{G}}_{c}(\th{x}_{M(c),q}; \mu_c),  \\
&\Bigl[{v}^{\bb}_{\rb,c} + v^{\gamma}_{\rb,c} \Bigr] \circ {\mathcal{G}}_{c}(\th{x}_{M(c),q}; \mu_c); \th{x}_{M(c),q},\mu_c \biggr) = 0,
\end{aligned}
\end{equation}
so that the system-level HRBE solution $\wtilde{u}_\rb(\mu) \in \mathcal{V}_\rb$ is given by $\wtilde{u}_\rb(\mu) = \sum_{c \in \sys} \left[ \wtilde{u}^\bb_{\rb,c}(\mu) + \wtilde{u}^{\gamma}_{\rb,c}(\mu) \right]$.

\section{Online-adaptive refinement through hierarchical error estimation}
\label{sec:adaptive}
Thus far, we have assumed that each archetype component in the library has a single RB space. As discussed in the introduction, a key challenge in component-based RB methods is determining the optimal fidelity (i.e., dimension) of each component/local RB model to achieve the desired system-level/global solution accuracy during the online phase. To address this challenge, we now assume that each archetype component has a set of RB spaces with varying fidelities and introduce a strategy to adaptively determine the appropriate fidelity of the local RB models for the instantiated components in the system so that the HRBE solution satisfies a prescribed error tolerance relative to the truth solution for the minimal computational cost. The proposed strategy follows a hierarchical error estimation framework, where solution accuracy is adaptively improved by selectively increasing the fidelities of the bubble and port RB spaces of the instantiated components.

\subsection{Multi-fidelity components, ports, and system}
\label{subsec:hierarchicallib}
We begin by introducing hierarchical RB spaces for the archetype ports. We recall the $N_\wp$-dimensional archetype port RB space $\what{\mathcal{X}}_{\rb,\wp}$ $\forall \wp \in \wpset$, introduced in Section~\ref{subsec:hrbe}. We refer to these as the \emph{finest} (i.e., the highest fidelity) RB spaces for the archetype ports. For each archetype port $\wp$, we introduce $N_{\wp}^{\mathrm{fdl}}$ hierarchical RB spaces $\what{\mathcal{X}}_{\rb,\wp}^1 \subset \cdots \subset \what{\mathcal{X}}_{\rb,\wp}^{N_{\wp}^{\mathrm{fdl}}} \equiv \what{\mathcal{X}}_{\rb,\wp}$ such that $\what{\mathcal{X}}_{\rb,\wp}^{f} = \text{span} \{ \what{\chi}_{\wp,i} \}_{i=1}^{N_{\wp,f}}$, $f \in \{ 1,\dots,N_{\wp}^{\mathrm{fdl}} \}$, where $N_{\wp,1} < \cdots < N_{\wp, N_{\wp}^{\mathrm{fdl}}} \equiv N_\wp$; i.e., we simply select a subset of the RB for the finest space to construct a \emph{coarse} space. We describe the procedure for constructing the finest space that induces this hierarchical structure in Section~\ref{subsec:rbtraining}.

We now introduce hierarchical RB spaces for the bubble RB space of the archetype components in an analogous manner. For each archetype component $\wc \in \wcset$, we recall the $N_\wc^\bb$-dimensional bubble RB space $\what{\mathcal{V}}_{\rb,\wc}^\bb$, introduced in Section~\ref{subsec:hrbe}, and refer to it as the finest bubble RB space for $\wc$. For each $\wc$, we introduce $N_{\wc, \bb}^{\mathrm{fdl}}$ hierarchical bubble RB spaces $\what{\mathcal{V}}_{\rb,\wc}^{\bb,1} \subset \cdots \subset \what{\mathcal{V}}_{\rb,\wc}^{\bb,N_{\wc, \bb}^{\mathrm{fdl}}} \equiv \what{\mathcal{V}}_{\rb,\wc}^{\bb}$ such that $\what{\mathcal{V}}_{\rb,\wc}^{\bb,f} = \text{span} \{ \what{\xi}_{\wc,i}^\bb \}_{i=1}^{N_{\wc,f}^{\bb}}$, $f \in \{ 1,\dots,N_{\wc, \bb}^{\mathrm{fdl}} \}$, where $N_{\wc,1}^{\bb} < \cdots < N_{\wc,N_{\wc, \bb}^{\mathrm{fdl}}}^{\bb} \equiv N_{\wc}^{\bb}$; i.e., analogously to the port spaces, we simply select a subset of the RB for the finest space to construct a coarse space. We describe the procedure for constructing the finest space that induces this hierarchical structure in Section~\ref{subsec:rbtraining}.

Subsequently, for each archetype component $\wc \in \wcset$, we introduce hierarchical port-lifted RB spaces $\what{\mathcal{V}}_{\rb,\wc}^{p,1} \subset \cdots \subset \what{\mathcal{V}}_{\rb,\wc}^{p,N_{\pi_{\wc}(p)}^{\mathrm{fdl}}} \equiv \what{\mathcal{V}}_{\rb,\wc}^{p}$, $\forall p \in \pset_\wc$, such that $\what{\mathcal{V}}_{\rb,\wc}^{p,f} = \text{span} \{ \what{\theta}_{\wc,i}^p \}_{i=1}^{N_{\pi_{\wc}(p),f}}$, $f \in \{ 1,\dots,N_{\pi_{\wc}(p)}^{\mathrm{fdl}} \}$. Additionally, for each $\wc$, we introduce $\{ \what{\mathcal{V}}_{\rb,\wc}^{\mathbf{f}} \}_{\mathbf{f} \in \mathcal{F}_\wc}$ as a family of multi-indexed RB spaces formed by employing one of the hierarchical bubble RB spaces and independently employing a combination of hierarchical port-lifted RB spaces across the component's ports. The index set $\mathcal{F}_\wc$ consists of all $(1 + n_\wc^\gamma)$-tuples $\mathbf{f} = (f_\bb, ( f_p)_{p \in \pset_\wc})$, where $f_\bb \in \{1,\dots,N_{\wc,\bb}^{\mathrm{fdl}} \}$ and $f_p \in \{1,\dots,N_{\pi_{\wc}(p)}^{\mathrm{fdl}} \}$ denote the fidelity level of the incorporated bubble and port $p$'s lifted RB spaces, respectively. We denote $N_\wc^{\mathrm{fdl}} \equiv N_{\wc,\bb}^{\mathrm{fdl}} \prod_{p \in \pset_\wc} N_{\pi_{\wc}(p)}^{\mathrm{fdl}}$ as the total cardinality of $\mathcal{F}_\wc$. Therefore, the previously introduced RB space $\what{\mathcal{V}}_{\rb,\wc} = \what{\mathcal{V}}_{\rb,\wc}^{(N_{\wc,\bb}^{\mathrm{fdl}}, (N_{\pi_{\wc}(p)}^{\mathrm{fdl}})_{p \in \pset_\wc} )}$ serves as the finest RB space for $\wc$, whereas the coarsest (i.e., the lowest fidelity) RB space for $\wc$ is $\what{\mathcal{V}}_{\rb,\wc}^{(1, \dots, 1 )}$. Each $\what{\mathcal{V}}_{\rb,\wc}^{\mathbf{f}}$ is associated with a RQ rule $(\th{x}_{\wc,q}^{\mathbf{f}}, \th{\rho}_{\wc,q}^{\mathbf{f}} )_{q=1}^{\wtilde{Q}_{\wc,\mathbf{f}}} \subset (\widehat{x}_{\wc,q}, \wrho_{\wc,q} )_{q=1}^{Q_\wc}$, obtained through the component-wise EQP outlined in Section~\ref{subsec:rqtraining}.

For each instantiated component $c \in \sys$, we define hierarchical RB spaces for their bubble and port spaces as
\begin{equation}\nonumber
  \begin{aligned}
    \mathcal{V}_{\rb,c}^{\bb,f} &\equiv \left\{v = \widehat{v} \circ \mathcal{G}_c^{-1}(\cdot; \mu_c) : \: \widehat{v} \in \widehat{\mathcal{V}}_{\rb,M(c)}^{\bb,f} \right\} \quadd f \in \{1, \dots, N_{M(c),\bb}^{\mathrm{fdl}} \}, \\
    \mathcal{V}_{\rb,c}^{p,f} &\equiv \left\{v = \widehat{v} \circ \mathcal{G}_c^{-1}(\cdot; \mu_c) : \: \widehat{v} \in \widehat{\mathcal{V}}_{\rb,M(c)}^{p, f} \right\} \quadd f \in \{1, \dots, N_{\pi_{M(c)}(p)}^{\mathrm{fdl}} \} \quad \forall p \in \pset_{M(c)}.
  \end{aligned}
\end{equation}
Additionally, for each $c \in \sys$, we introduce the set of $N_{M(c)}^{\mathrm{fdl}}$ multi-indexed RB spaces $\mathcal{V}_{\rb,c}^{\mathbf{f}} \equiv \Big\{v = \widehat{v} \circ \mathcal{G}_c^{-1}(\cdot; \mu_c) : \: \widehat{v} \in \widehat{\mathcal{V}}_{\rb,M(c)}^{\mathbf{f}} \Big\}$, $\mathbf{f} \in \mathcal{F}_{M(c)}$. Subsequently, the system's finest RB space, denoted by $\mathcal{V}_\rb^{\mathrm{finest}}$, is constructed by assembling the finest RB spaces of all the instantiated components, whereas the coarsest RB space, denoted by $\mathcal{V}_\rb^{\mathrm{coarsest}} \subset \mathcal{V}_\rb^{\mathrm{finest}}$, is formed by assembling the coarsest RB spaces of all the instantiated components.

\subsection{Hierarchical error estimation}
We now introduce a hierarchical error estimator to approximate the system-level error in the HRBE solution relative to the truth solution and identify the components that require higher fidelity RB spaces to improve the solution accuracy. 

For any $\wc \in \wcset$, we introduce $\bar{\mathcal{F}}_\wc \equiv \mathcal{F}_\wc \setminus \{(N_{\wc, \bb}^{\mathrm{fdl}}, (N_{\pi_\wc(p)}^{\mathrm{fdl}})_{p \in \pset_{\wc}})\}$; i.e., a set of all multi-indexed fidelities except the finest. We assume at the $k$-th adaptive refinement iteration, $\mathcal{V}_{\rb}^k \equiv \left( \boplus_{c \in \sys} \mathcal{V}_{\rb,c}^{\mathbf{f}_c} \right) \cap \mathcal{V}$ denotes the intermediate RB space for the system, where $\mathbf{f}_c = (f_{c, \bb}, (f_{c, p})_{p \in \pset_{M(c)}}) \in \bar{\mathcal{F}}_{M(c)}$ and $\mathcal{V}_\rb^{\mathrm{coarsest}} \subset \mathcal{V}_\rb^k \subset \mathcal{V}_\rb^{\mathrm{finest}}$. At this iteration, we introduce the \emph{refined} system-level RB space ${\mathcal{V}}_{\rb}^{k'} = \left( \boplus_{c \in \sys} \mathcal{V}_{\rb,c}^{\mathbf{f}_c'} \right) \cap \mathcal{V}$, where $\mathbf{f}_c' = \mathbf{f}_c + \mathbf{1} \equiv (f_{c, \bb} + 1, (f_{c, p} + 1)_{p \in \pset_{M(c)}})$ $\forall c \in \sys$ so that all bubble and port fidelities are incremented by one. We now introduce the following proposition, which formulates a system-level error estimator by extending the hierarchical error estimator proposed in~\cite{hain2019hierarchical} to component-based systems.

\begin{proposition}\label{prop:errestimator} For a given $\mu \in \mathcal{D}$, at the $k$-th adaptive refinement iteration, let $\wtilde{u}_\rb(\mu) \in \mathcal{V}_\rb^k$ and $\wtilde{{u}}_\rb'(\mu) \in {\mathcal{V}}_\rb^{k'}$ be the HRBE solutions obtained by solving~\eqref{eq:nlhrbepr} using the RB spaces $\mathcal{V}_\rb = \mathcal{V}_{\rb}^k$ and $\mathcal{V}_\rb =  {\mathcal{V}}_\rb^{k'}$, respectively. We introduce non-negative constants $\eta_{\wc,\mathbf{f}}$ for each $\wc \in \wcset$ and $\mathbf{f} \in \bar{\mathcal{F}}_{\wc}$ such that
  \aequation{\label{eq:eta}
    \|u_h(\mu)\big|_{\Omega_c} - \wtilde{{u}}_\rb'(\mu)\big|_{\Omega_c} \|_{\mathcal{V}_c} \leq \eta_{M(c), {\mathbf{f}}_c} \: \|u_h(\mu)\big|_{\Omega_c} - \wtilde{u}_\rb(\mu)\big|_{\Omega_c} \|_{\mathcal{V}_c}.
  }
  Since $\mathcal{V}_{\rb}^k \subset {\mathcal{V}}_\rb^{k'}$, we assume that $\wtilde{{u}}_\rb'(\mu)$ provides a more accurate approximation of the truth solution $u_h(\mu) \in \mathcal{V}_h$ than $\wtilde{u}_\rb(\mu)$. Accordingly, we assume $0 \leq \eta_{\wc,\mathbf{f}} < 1$ for all $\wc \in \wcset$ and $\mathbf{f} \in \bar{\mathcal{F}}_\wc$. The system-level error can then be estimated as  
  \aequation{\label{eq:errestimator}
    \|u_h(\mu) - \wtilde{u}_\rb(\mu)\|_{\mathcal{V}} \leq \sum_{c \in \sys} \frac{1}{1 - \eta_{M(c), {\mathbf{f}}_c}} \|\wtilde{{u}}_\rb'(\mu)\big|_{\Omega_c} - \wtilde{u}_\rb(\mu)\big|_{\Omega_c} \|_{\mathcal{V}_c}.
  }
\end{proposition}
\begin{proof}
  Under the assumption that $0 \leq \eta_{\wc,\mathbf{f}} < 1$ $\forall \wc \in \wcset$ and $\forall \mathbf{f} \in \bar{\mathcal{F}}_\wc$, for any $c \in \sys$, we have
  \aequation{\label{eq:errestimator1}
    (1 - \eta_{M(c), {\mathbf{f}}_c}) \|u_h(\mu)\big|_{\Omega_c} - \wtilde{u}_\rb(\mu)\big|_{\Omega_c} \|_{\mathcal{V}_c} &\leq \|{u}_h(\mu)\big|_{\Omega_c} - \wtilde{u}_\rb(\mu)\big|_{\Omega_c} \|_{\mathcal{V}_c} - \|{u}_h(\mu)\big|_{\Omega_c} - \wtilde{{u}}_\rb'(\mu)\big|_{\Omega_c} \|_{\mathcal{V}_c}
    \\&\leq \|\wtilde{{u}}_\rb'(\mu)\big|_{\Omega_c} - \wtilde{u}_\rb(\mu)\big|_{\Omega_c} \|_{\mathcal{V}_c},
  }
  where the first inequality follows from subtracting $\|u_h(\mu)\big|_{\Omega_c} - \wtilde{u}_\rb(\mu)\big|_{\Omega_c} \|_{\mathcal{V}_c}$ from both sides of~\eqref{eq:eta} and multiplying by $-1$, and the second inequality follows from the triangular inequality. Therefore, incorporating~\eqref{eq:errestimator1} leads to
  \begin{equation}\nonumber
    \begin{aligned}
      \|u_h(\mu) - \wtilde{u}_\rb(\mu)\|_{\mathcal{V}} &\leq \sum_{c \in \sys} \|u_h(\mu)\big|_{\Omega_c} - \wtilde{u}_\rb(\mu)\big|_{\Omega_c} \|_{\mathcal{V}_c} \leq \sum_{c \in \sys} \frac{1}{1 - \eta_{M(c), {\mathbf{f}}_c}} \|\wtilde{{u}}_\rb'(\mu)\big|_{\Omega_c} - \wtilde{u}_\rb(\mu)\big|_{\Omega_c} \|_{\mathcal{V}_c},
    \end{aligned}
  \end{equation}
  where the first inequality follows from the triangular inequality.
\end{proof}

We refer to $\eta_{\wc, \mathbf{f}}$ $\forall \wc \in \wcset$ and $\forall \mathbf{f} \in \bar{\mathcal{F}}_\wc$ in Proposition~\ref{prop:errestimator} as the \emph{error contraction factors}. The error contraction factor $\eta_{\wc, \mathbf{f}}$ for $\wc$ quantifies the relative reduction in its HRBE solution error when its RB space is refined from $\mathcal{V}_{\rb, \wc}^{\mathbf{f}}$ to $\mathcal{V}_{\rb, \wc}^{\mathbf{f}'}$, where $\mathbf{f}' = \mathbf{f} + \mathbf{1}$. Smaller values of $\eta_{\wc, \mathbf{f}}$ indicate a more significant improvement in accuracy due to refinement, whereas values closer to one suggest diminishing returns from further enrichment. The computational procedure for determining error contraction factors during the offline phase is described in Section~\ref{subsec:errorcontraction}; for now, we assume the error contraction factors are known.

In addition to the (global) error estimate for the assembled system, Proposition~\ref{prop:errestimator} provides component-wise local error indicators that inform the adaptive refinement of the bubble and port-lifted RB spaces for individual components within the system. At each refinement iteration, by evaluating the error between the current HRBE solution and its refined counterpart, we can estimate the error in the former and identify components that contribute most significantly to the system-level error. This localized error information guides the refinement process by prioritizing the components that require further enrichment. This process is analogous to adaptive finite element methods guided by element-wise error indicators.

\subsection{Adaptive refinement strategy}
We now present a strategy for adaptively refining the component RB spaces during the online phase to obtain an HRBE solution that satisfies the prescribed solution error relative to the truth solution. The algorithm begins with the coarsest possible RB space, where each component is initialized with its lowest-fidelity bubble and port-lifted RB spaces, i.e., ${\mathbf{f}}_c = (f_{c, \bb}, (f_{c, p})_{p \in \pset_{M(c)}}) = (1, \dots, 1)$ $\forall c \in \sys$. The refinement process is guided by a system-level error estimator based on Proposition~\ref{prop:errestimator}, which quantifies the contribution of each component to the overall solution error.

At the $k$-th iteration, the current HRBE solution $\wtilde{u}_\rb(\mu) \in \mathcal{V}_\rb^k$ is computed by solving the HRBE problem~\eqref{eq:nlhrbepr} using $\mathcal{V}_\rb = \mathcal{V}_\rb^k = \left( \boplus_{c \in \sys} \mathcal{V}_{\rb,c}^{\mathbf{f}_c} \right) \cap \mathcal{V}$, $\mathbf{f}_c \in \bar{\mathcal{F}}_{M(c)}$. To estimate the relative error in this solution, the algorithm constructs the refined RB space ${\mathcal{V}}_\rb^{k'} = \left( \boplus_{c \in \sys} \mathcal{V}_{\rb,c}^{\mathbf{f}_c'} \right) \cap \mathcal{V}$, where $\mathbf{f}_c' = \mathbf{f}_c + \mathbf{1}$ $\forall c \in \sys$. The refined HRBE solution $\wtilde{{u}}_\rb'(\mu) \in {\mathcal{V}}_\rb^{k'}$ is then computed using ${\mathcal{V}}_\rb = {\mathcal{V}}_\rb^{k'}$ in~\eqref{eq:nlhrbepr}. The system-level error estimate $\mathcal{E}$ is computed using~\eqref{eq:errestimator}, the error contraction factors $\eta_{M(c), \mathbf{f}_c}$, and $\|\wtilde{{u}}_\rb'(\mu)\big|_{\Omega_c} - \wtilde{u}_\rb(\mu)\big|_{\Omega_c}\|_{\mathcal{V}_c}$ $\forall c \in \sys$. As in the case study presented in Section~\ref{sec:res}, we may approximate the relative error in the current HRBE solution by $\mathcal{E} / \| \wtilde{{u}}_\rb'(\mu) \|_{\mathcal{V}}$, where $\| \wtilde{{u}}_\rb'(\mu) \|_{\mathcal{V}}$ is used in place of $\|u_h(\mu)\|_{\mathcal{V}}$. These two quantities are expected to converge as the refinement proceeds.

If the estimated error satisfies the prescribed tolerance, the procedure terminates and the current HRBE solution is returned. Otherwise, the algorithm identifies the subset of components that contribute most significantly to the error. Specifically, the components are ranked based on their local error contributions, and the top $\Delta\%$ are selected for refinement. For each selected component, the bubble and port-lifted RB spaces are increased in fidelity by one level, provided they have not yet reached their \emph{allowable} finest RB space. The refinement process continues for a maximum of $N_{\mathrm{ref}}$ iterations or until convergence is achieved. Algorithm~\ref{alg:adaptiveref} summarizes the described adaptive refinement strategy.

\begin{remark}\label{rem:refinement}
  In the hierarchical error estimation framework adopted here, the error in the current RB solution is estimated by comparing it against a solution computed with higher-fidelity bubble and port-lifted RB spaces for all components. As a result, for all components, we reserve the finest available bubble and port-lifted RB spaces to serve as the comparison space during refinement. Therefore, during refinement, the allowable fidelity level of the bubble and port-lifted RB spaces $\forall c \in \sys$ and $\forall p \in \pset_{M(c)}$ is effectively $N_{M(c), \bb}^{\mathrm{fdl}} - 1$ and $N_{\pi_{M(c)}(p)}^{\mathrm{fdl}} - 1$, respectively. This restriction ensures that finer bubble and port-lifted RB spaces always exist to enable error estimation.
\end{remark}

\begin{algorithm}[!tb]
  \caption{Adaptive refinement of component RB spaces during the online phase to find an HRBE solution that satisfies the prescribed accuracy.}\label{alg:adaptiveref}
  \SetAlgoLined
  \KwInput{System-level parameter $\mu \in \mathcal{D}$; desired HRBE solution error $\epsilon > 0$; maximum number of refinement iterations $N_{\mathrm{ref}}$; percentage of components to refine per iteration $\Delta$}
  \KwOutput{If converged, an HRBE solution $\wtilde{u}_\rb(\mu)$ that satisfies $\|u_h(\mu) -  \wtilde{u}_\rb(\mu)\|_{\mathcal{V}} \leq \epsilon$}
  Initialize system-level error estimate: $\mathcal{E} = \infty$\;
  Initialize refinement iteration counter: $k = 0$\;
  Initialize the current fidelity levels: ${\mathbf{f}}_c = (f_{c, \bb}, (f_{c, p})_{p \in \pset_{M(c)}}) = (1, \dots, 1)$ $\forall c \in \sys$\;
  Initialize $\sys_{\mathrm{sub}} = \emptyset$ \tcp*{Subset of components selected for refinement}
  \While{$k < N_{\mathrm{ref}}$}{ 
      \tcp{Solve the HRBE problem with current RB spaces}
      Solve the HRBE problem~\eqref{eq:nlhrbepr} using $\mathcal{V}_\rb = \mathcal{V}_\rb^k = \left( \boplus_{c \in \sys} \mathcal{V}_{\rb,c}^{\mathbf{f}_c} \right) \cap \mathcal{V}$ to obtain $\wtilde{u}_\rb(\mu) \in \mathcal{V}_\rb^k$\;
      \tcp{Check if further refinement is possible}
      \If{$f_{c,\bb} = N_{M(c), \bb}^{\mathrm{fdl}} - 1$ $\forall c \in \sys$ and $f_{c,p} = N_{\pi_{M(c)}(p)}^{\mathrm{fdl}} - 1$ $\forall p \in \pset_{M(c)}$}{
          \textbf{break}\;
      }
      \tcp{Solve the refined HRBE problem}
      $\mathbf{f}_c' = \mathbf{f}_c + \mathbf{1}$ $\forall c \in \sys$\;
      Solve the HRBE problem using ${\mathcal{V}}_\rb = {\mathcal{V}}_\rb^{k'} = \left( \boplus_{c \in \sys} \mathcal{V}_{\rb,c}^{\mathbf{f}_c'} \right) \cap \mathcal{V}$ to obtain $\wtilde{{u}}_\rb'(\mu) \in {\mathcal{V}}_\rb^{k'}$\;

      \tcp{Estimate the system-level error}
      Compute component-wise errors $\mathcal{E}_c \equiv \|\wtilde{{u}}_\rb'(\mu)\big|_{\Omega_c} - \wtilde{u}_\rb(\mu)\big|_{\Omega_c}\|_{\mathcal{V}_c} / (1 - \eta_{M(c), {\mathbf{f}}_c})$ $\forall c \in \sys$\;
      Compute the system-level error estimate $\mathcal{E} = \sum_{c \in \sys} \mathcal{E}_c$\;
      \If{$\mathcal{E} \leq \epsilon$}{
          \textbf{break}\;
      }
      \tcp{Select components for refinement}
      Sort $\mathcal{E}_c$ $\forall c \in \sys$ in descending order\;
      Identify the top $\Delta\%$ of components with the highest local errors for refinement to form $\sys_{\mathrm{sub}} \subset \sys$\;
      \tcp{Refine the RB spaces}
      $\mathbf{f}_c \leftarrow (\max(f_{c, \bb} + 1, N_{M(c),\bb}^{\mathrm{fdl}} - 1), (\max(f_{c, p} + 1, N_{\pi_{M(c)(p)}}^{\mathrm{fdl}} - 1))_{p \in \pset_{M(c)}})$ $\forall c \in \sys_{\mathrm{sub}}$\;
      $k \leftarrow k + 1$\;
  }    
\end{algorithm}

\subsection{Computational cost and memory requirement}
\label{subsec:costmemory}
We first compare the computational cost of solving the truth problem with that of solving the HRBE problem using the adaptive refinement strategy in Algorithm~\ref{alg:adaptiveref}. We use Newton's method to solve~\eqref{eq:truthref} and~\eqref{eq:nlhrbepr}. Each Newton step for solving~\eqref{eq:truthref} involves~$\mathcal{O}({Q}_h \equiv \sum_{c \in \sys} Q_{M(c)})$ operations to evaluate the truth residual and Jacobian. It also requires solving a sparse linear system of equations, which entails~$\mathcal{O}(\mathcal{N}_h^n)$ operations, where $1 \leq n \leq 2$ depends on the domain dimension and the choice of solver. On the other hand, the $k$-th adaptive refinement iteration involves solving the HRBE problem~\eqref{eq:nlhrbepr} using the current RB space~$\mathcal{V}_\rb^k$ and the refined RB space~${\mathcal{V}}_\rb^{k'}$. In each Newton step, evaluating the residual and Jacobian using $\mathcal{V}_\rb^k$ requires $\mathcal{O}(\sum_{c \in \sys} N_{c,k}^2 \wtilde{Q}_{c,k}) \ll \mathcal{O}(Q_h)$, where $N_{c,k}$ denotes the size of the current RB space for component $c$, and $\wtilde{Q}_{c,k}$ denotes the number of RQ points for component $c$ at this iteration. Evaluating the residual and Jacobian using ${\mathcal{V}}_\rb^{k'}$ proceeds similarly. Additionally, solving the linear system at each Newton step using $\mathcal{V}_\rb^k$ and ${\mathcal{V}}_\rb^{k'}$ requires $\mathcal{O}([\dim(\mathcal{V}_\rb^k)]^m) \leq \mathcal{O}([\dim({\mathcal{V}}_\rb^{k'})]^m) \ll \mathcal{O}(\mathcal{N}_h^n)$ operations, where $1 \leq m \leq 3$ depends on the domain dimension, number of components, and the choice of solver. Specifically, $1 \leq m \leq 2$ when the number of components $N_\comp$ is large and thus the Jacobian is component-block-wise sparse; in the extreme case where $N_\comp = 1$ and the Jacobian is dense, $m = 3$. The remaining operations in Algorithm~\ref{alg:adaptiveref} incur negligible computational cost compared to solving the two HRBE problems.

We now comment on the memory requirements. The memory footprint of the truth problem~\eqref{eq:truthref}, dominated by the storage of the truth Jacobian, is $\mathcal{O}\left( \mathcal{N}_h^n \right)$, where $n = 1$ if an iterative linear solver is used at each Newton step, and $n = 4/3$ in the worst case for $d = 3$ for storing the factorization using a sparse direct solver. For the online phase, the following data must be loaded into memory for each archetype component $\wc \in \wcset$: (i) the error contraction factor $\eta_{\wc, \mathbf{f}}$ $\forall \mathbf{f} \in \bar{\mathcal{F}}_\wc$ (ii) the RQ weights $\{ \th{\rho}_{\wc,q}^{\mathbf{f}}\}_{q=1}^{\wtilde{Q}_{\wc,\mathbf{f}}}$ $\forall \mathbf{f} = (f_\bb, ( f_p)_{p \in \pset_\wc}) \in {\mathcal{F}}_\wc$, (iii) the value of the bubble basis functions $\{ \widehat{\xi}_{\wc,i}^\bb \}_{i=1}^{N_{\wc,f_{\bb}}^\bb}$ and port-lifted basis functions $\{ \widehat{\theta}_{\wc,i}^p \}_{i=1}^{N_{\pi_\wc(p),f_p}}$ $\forall p \in \pset_\wc$ at the RQ points $\{\th{x}_{\wc,q}^{\mathbf{f}} \}_{q=1}^{\wtilde{Q}_{\wc,\mathbf{f}}}$ $\forall \mathbf{f}$, and (iv) the gradients of the bubble basis functions $\{ \nabla \widehat{\xi}_{\wc,i}^\bb \}_{i=1}^{N_{\wc,f_\bb}^\bb}$ and port-lifted basis functions $\{ \nabla \widehat{\theta}_{\wc,i}^p \}_{i=1}^{N_{\pi_\wc(p),f_p}}$ $\forall p \in \pset_\wc$, also evaluated at the RQ points $\{\th{x}_{\wc,q}^{\mathbf{f}} \}_{q=1}^{\wtilde{Q}_{\wc,\mathbf{f}}}$ $\forall \mathbf{f} \in {\mathcal{F}}_\wc$. Therefore, the total memory footprint for loading the library data is
\aequation{\nonumber
  \mathcal{M}_{\mathrm{lib}} \equiv
   \sum_{\wc \in \wcset} \Big( (N_\wc^{\mathrm{fdl}} - 1) + \sum_{\mathbf{f} \in \mathcal{F}_\wc} \wtilde{Q}_{\wc, \mathbf{f}} \Big(1 + (N_{\wc, f_\bb}^\bb + \sum_{p \in \pset_\wc} N_{\pi_\wc(p),f_p})(d + 1) \Big) \Big).
}
Additionally, the storage of the system-level HRBE residual and Jacobian at the $k$-th refinement iteration requires $\mathcal{O}(\sum_{c \in \sys} N_{c,k'}^2)$ memory, where $N_{c,k'}$ denotes the size of the refined RB space for component $c$ at this iteration. The total memory footprint of the online phase is thus $\mathcal{O}\left(\mathcal{M}_{\mathrm{lib}} + \max_{k' \in \{1, \dots, N_{\mathrm{ref}} \}} \sum_{c \in \sys} N_{c,k'}^2 \right)$, which is independent of $\mathcal{N}_\wc^\bb$ $\forall \wc \in \wcset$, $\mathcal{N}_\wp$ $\forall \wp \in \wpset$, and $Q_\wc$ $\forall \wc \in \wcset$.

\section{Component-wise offline training}
\label{sec:training}
In this section, we describe the computational procedures required in the offline phase to support the adaptive refinement strategy introduced in Section~\ref{sec:adaptive}. We begin by outlining the construction of the RB spaces with varying fidelities for the archetype components and ports. We then describe the component-wise EQP for determining the RQ rules associated with each RB space fidelity. Finally, we present a methodology for estimating the error contraction factors.

\subsection{Component-wise RB construction}
\label{subsec:rbtraining}
We first introduce a routine for generating truth snapshot solutions for the archetype components and ports in the library. We follow the procedure introduced in~\cite{ebrahimi2024hyperreduced}, which belongs to the family of subsystem-based training procedures for component-based systems, such as the pairwise training procedure for port modes in~\cite{eftang2013port, eftang2014port}. We note that the procedure in~\cite{ebrahimi2024hyperreduced} for nonlinear systems considers a larger subsystem than the procedures for linear systems, since the training set needs to cover not only relevant solution \emph{shapes/modes} but also representative solution \emph{magnitudes}. 

For each archetype component $\wc \in \wcset$, we define a parameter training set $\Xi^{\train}_\wc \equiv \{ \mu^{\train}_{\wc,n} \in \widehat{\mathcal{D}}_\wc \}_{n=1}^{N^{\train}_\wc}$, where $N^{\train}_\wc$ denotes the number of training parameter samples. For each $\wc \in \wcset$, we generate $N_{\sample}$ sample \emph{subsystems} by connecting it to other randomly selected components from the library via its $n^\gamma_\wc$ local ports. The probability of establishing a connection through each local port is given by $\nu$. We assign random parameter values to each component in the assembled subsystems from their corresponding parameter training sets and apply independent random constant Dirichlet boundary conditions, with uniform density, to all boundary global ports. We then solve the truth problem~\eqref{eq:truthref} for the subsystems and extract the truth solutions on $\wc$. The extracted solutions are added to the truth snapshot set $U_{h,\wc}^{\train}$ associated with component $\wc$.

We then decompose the extracted solutions into their bubble and port parts as in~\eqref{eq:compsol}. The bubble parts are mapped to the reference domain and added to the bubble snapshot set of the corresponding archetype component. Similarly, the port parts are mapped and added to the snapshot sets of their corresponding archetype ports. We introduce $U_{h,\wc}^{\train,\bb}$ $\forall \wc \in \wcset$ and $U_{h,\wp}^{\train}$ $\forall \wp \in \wpset$ as the snapshot sets for the archetype components' bubble and archetype ports, respectively. The core assumption of this procedure, outlined in Algorithm~\ref{alg:emptr}, is that the generated snapshot solutions adequately capture the range of possible solutions that components and ports may encounter in actual system configurations. 

Once the snapshot sets are formed, we construct hierarchical RB spaces using proper orthogonal decomposition (POD) with varying levels of accuracy. For each archetype component $\wc \in \wcset$, we apply POD to $U_{h,\wc}^{\train,\bb}$ with decreasing tolerances $\delta_{\wc,\mathrm{pod}}^1 > \cdots > \delta_{\wc,\mathrm{pod}}^{N_{\wc,\bb}^{\mathrm{fdl}}}$ to find hierarchical bubble RB spaces $\{\what{\mathcal{V}}_{\rb,\wc}^{\bb,f}\}_{f=1}^{N_{\wc,\bb}^{\mathrm{fdl}}}$. Similarly, for each archetype port $\wp \in \wpset$, we apply POD to $U_{h,\wp}^\train$ with decreasing tolerances $\delta_{\wp,\mathrm{pod}}^1 > \cdots > \delta_{\wp,\mathrm{pod}}^{N_\wp^{\mathrm{fdl}}}$ to find hierarchical RB spaces $\{\what{\mathcal{X}}_{\rb,\wp}^f\}_{f=1}^{N_\wp^{\mathrm{fdl}}}$. Subsequently, for each $\wc \in \wcset$, we construct a hierarchy of port-lifted RB spaces $\{\what{\mathcal{V}}_{\rb,\wc}^{p,f}\}_{f=1}^{N_{\pi_\wc(p)}^{\mathrm{fdl}}}$ for each port $p \in \pset_\wc$ by elliptically lifting the port RB spaces using~\eqref{eq:portlift}. Finally, we assemble the bubble and port-lifted RB spaces to obtain the component RB spaces $\{\what{\mathcal{V}}_{\rb,\wc}^{\mathbf{f}}\}_{\mathbf{f} \in \mathcal{F}_\wc}$, as previously introduced in Section~\ref{subsec:hierarchicallib}.

\begin{algorithm}[!tb]
  \caption{Generating snapshot solutions for component-wise training.}\label{alg:emptr}
  \SetAlgoLined
  \KwInput{Number of sample subsystems $N_{\sample}$; probability of port connection $0 \leq \nu \leq 1$}
  \KwOutput{Snapshot sets $U_{h,\wc}^{\train}$ and $U_{h,\wc}^{\train,\bb}$ $\forall \wc \in \wcset$ and $U_{h,\wp}^{\train}$ $\forall \wp \in \wpset$}
  Initialize $U_{h,\wc}^{\train} \leftarrow \emptyset$ and $U_{h,\wc}^{\train,\bb} \leftarrow \emptyset$ $\forall \wc \in \wcset$\;
  Initialize $U_{h,\wp}^{\train} \leftarrow \emptyset$ $\forall \wp \in \wpset$\;
  \ForEach{$\wc \in \wcset$}{
      \For{$n=1, \dots, N_{\sample}$}{
        \tcp{Assemble a sample subsystem $\mathcal{C}_{\rm sub}$ for $\wc$}
        Initialize $\mathcal{C}_{\rm sub} \leftarrow \{ \wc \}$\;
        \ForEach{$p \in \pset_\wc$}{
          With probability $\nu$, connect port $p$ of $\wc$ to another randomly selected component in the library\;
          Add the selected component to $\mathcal{C}_{\rm sub}$\;
        }
        \tcp{Assign parameters and boundary conditions}
        Assign $\mu_c \sim \text{Uniform}(\mathcal{D}_c)$ to each $c \in \mathcal{C}_{\rm sub}$\;
        Apply random constant Dirichlet boundary conditions to all boundary global ports\;
        \tcp{Solve the subsystem and extract the solution for $\wc$}
        Solve the truth problem~\eqref{eq:truthref} on $\mathcal{C}_{\rm sub}$\;
        Extract the truth solution $u_{h,c}^{\train}$ on component $c$ that is instantiated from $\wc$\;
        Set $u_{h,\wc}^{\train} = u_{h,c}^{\train} \circ \mathcal{G}_c(\cdot; \mu_c)$\;
        $U_{h,\wc}^{\train} \leftarrow U_{h,\wc}^{\train} \cup \{ u_{h,\wc}^{\train} \}$\;
        Decompose: $u_{h,c}^{\train} = u_{h,c}^{\train,\bb} + \sum_{p \in \pset_{M(c)}} u_{h,c}^{\train,p}$\;
        \tcp{Map and store the bubble part}
        Set $u_{h,\wc}^{\train,\bb} = u_{h,c}^{\train,\bb} \circ \mathcal{G}_c(\cdot; \mu_c)$\;
        $U_{h,\wc}^{\train,\bb} \leftarrow U_{h,\wc}^{\train} \cup \{ u_{h,\wc}^{\train,\bb} \}$\;
        \tcp{Map and store the port parts}
        \ForEach{$p \in \pset_\wc$}{
          Set $u_{h,\wc}^{\train,p} = u_{h,c}^{\train,p} \circ \mathcal{G}_c(\cdot; \mu_c)$\;
          Set $\bar{u}_{h,\wc}^{\train,p} = u_{h,\wc}^{\train,p} \big|_{\what{\gamma}_{\wc,p}} \circ \mathcal{R}_{\wc,p}(\cdot)$\;
          $U_{h,\pi_\wc(p)}^{\train} \leftarrow U_{h,\pi_\wc(p)}^{\train} \cup \{ \bar{u}_{h,\wc}^{\train,p} \}$\;
        }
      }
  }
\end{algorithm}

\subsection{Component-wise hyperreduction}
\label{subsec:rqtraining}
We employ the component-wise EQP~\cite{ebrahimi2024hyperreduced} to construct the RQ rules for each archetype component; see~\cite{ebrahimi2024hyperreduced} for a detailed description of the method. The EQP includes an accuracy hyperparameter for each archetype component $\wc \in \wcset$, which controls the accuracy of the component residual and Jacobian integrals computed with the RQ rule relative to their computation with the truth quadrature rule. In~\cite{ebrahimi2024hyperreduced}, we propose an offline training procedure that performs hyperreduction at multiple accuracy levels for all components. During the online phase, the appropriate hyperreduction accuracy and the corresponding RQ rule are selected adaptively to ensure that the HRBE solution satisfies a user-prescribed system-level error tolerance relative to the RB solution.

In this work, however, our goal is to control the error between the HRBE and truth solutions, which comprises both the RB approximation error and the hyperreduction error. Since beyond a certain hyperreduction accuracy the dominant source of error between HRBE and truth solutions is due to the RB approximation, we postulate that it is sufficient to ensure that the hyperreduction error remains below the RB approximation error. Therefore, we propose a strategy to couple hyperreduction and POD accuracies (i.e., fidelities) in the component-wise offline training. In the online phase, as described in Section~\ref{sec:adaptive}, given the desired system-level error, we adaptively select the RB fidelity of the components while allowing their hyperreduction fidelity to follow accordingly. This approach significantly reduces the computational complexity of the online phase by eliminating the need to adapt both RB and hyperreduction fidelities independently. 

\subsubsection{Error bound for the HRBE solution}
In this section, we derive a bound for the error between the HRBE and truth solutions. This bound is later used in Section~\ref{subsubsection:coupling} to develop a coupling strategy between hyperreduction and POD fidelities. For brevity and clarity, we defer the theoretical details to Section~\ref{sec:error}. Furthermore, to streamline the presentation, we temporarily disregard the multi-indexed RB space structure and assume that each component is associated with a single RB space, as introduced in Section~\ref{sec:hrbe}. The necessary extensions to account for the multi-indexed RB spaces are discussed in the following section.

We begin by introducing $\what{R}_{h,\wc}: \what{\mathcal{V}}_{h,\wc} \times \what{\mathcal{V}}_{h,\wc} \times \what{\mathcal{D}}_\wc \to \mathbb{R}$ and $\wtilde{\what{R}}_{\rb,\wc}: \what{\mathcal{V}}_{\rb,\wc} \times \what{\mathcal{V}}_{\rb,\wc} \times \what{\mathcal{D}}_\wc \to \mathbb{R}$ $\forall \wc \in \wcset$, respectively, as the reference-domain truth and HRBE residual forms such that
\aequation{\nonumber
  R_h(w, v;\mu) &= \sum_{c \in \sys} \what{R}_{h,M(c)}(w \big|_{\Omega_c} \circ \mathcal{G}_{c}(\cdot; \mu_c), v \big|_{\Omega_c} \circ \mathcal{G}_{c}(\cdot; \mu_c);\mu_c) \quadd \forall w,v \in \mathcal{V}_h, \forall \mu \in \mathcal{D},\\
  \wtilde{R}_\rb(w, v;\mu) &= \sum_{c \in \sys} \wtilde{\what{R}}_{\rb,M(c)}(w \big|_{\Omega_c} \circ \mathcal{G}_{c}(\cdot; \mu_c), v \big|_{\Omega_c} \circ \mathcal{G}_{c}(\cdot; \mu_c);\mu_c) \quadd \forall w,v \in \mathcal{V}_\rb, \forall \mu \in \mathcal{D}.
}
Explicit expressions for $\what{R}_{h,\wc}(\cdot, \cdot; \cdot)$ and $\wtilde{\what{R}}_{\rb,\wc}(\cdot, \cdot; \cdot)$ are provided in Appendix~\ref{app:derivative}. We assume that for each $\wc \in \wcset$, there exists some $\varepsilon_{\rb, \wc} \in \mathbb{R}_{\geq 0}$ and $\wtilde{\varepsilon}_{\rb, \wc} \in \mathbb{R}_{\geq 0}$ so that for the system $\sys$ and $\forall \mu = (\mu_c)_{c \in \sys} \in \mathcal{D}$, the following inequalities hold:
\begin{align}
  \sup_{v \in \mathcal{V}_{h,c}} \frac{|\what{R}_{h,M(c)}({u}_\rb(\mu)\big|_{\Omega_c} \circ \mathcal{G}_c(\cdot; \mu_c), {v} \circ \mathcal{G}_c(\cdot; \mu_c); \mu_c) |}{\| v \|_{\mathcal{V}_c}} &\leq \varepsilon_{\rb, M(c)} \quadd \forall c \in \sys \label{eq:epsilonhc}, \\
  \sup_{v \in \mathcal{V}_{\rb,c}} \frac{|\wtilde{\what{R}}_{\rb,M(c)}({u}_\rb(\mu)\big|_{\Omega_c} \circ \mathcal{G}_c(\cdot; \mu_c), {v} \circ \mathcal{G}_c(\cdot; \mu_c); \mu_c) |}{\| v \|_{\mathcal{V}_c}} &\leq \wtilde{\varepsilon}_{\rb, M(c)} \quadd \forall c \in \sys. \label{eq:epsilonrbc}
\end{align}
Then, under the regularity conditions elaborated in Section~\ref{sec:error}, the error between the truth and HRBE solutions can be bounded as
\aequation{\label{eq:totalerror}
      \|u_h(\mu) - \wtilde{u}_\rb(\mu) \|_{\mathcal{V}} \leq \frac{2}{\beta(\mu)} \sum_{c \in \sys} (\varepsilon_{\rb,M(c)}+ \wtilde{\varepsilon}_{\rb, M(c)}),
}
where $\beta(\mu) = \min (\beta_h({u}_\rb(\mu); \mu), \beta_\rb({u}_\rb(\mu); \mu))$ for the truth and RB inf--sup constants
\aequation{\nonumber
  \beta_h(v; \mu) = \inf_{w \in \mathcal{V}_h} \sup_{z \in \mathcal{V}_h} \frac{|R_h'(v, w, z; \mu)|}{\| w\|_{\mathcal{V}} \|z\|_{\mathcal{V}}} \quadd \forall v \in \mathcal{V}_h, \forall \mu \in \mathcal{D}, \\
  \beta_\rb(v; \mu) = \inf_{w \in \mathcal{V}_\rb} \sup_{z \in \mathcal{V}_\rb} \frac{|R_h'(v, w, z; \mu)|}{\| w\|_{\mathcal{V}} \|z\|_{\mathcal{V}}} \quadd \forall v \in \mathcal{V}_\rb, \forall \mu \in \mathcal{D}.
}
Here, $R'_h(v, w, z; \mu)$ $\forall v,w,z \in \mathcal{V}_h$ and $\forall \mu \in \mathcal{D}$ is the G\^{a}teaux derivative of the truth residual form $R_h(\cdot,z;\mu)$ at $v$ in the direction of $w$. Explicit expression for $R'_h(v,w,z; \mu)$ is provided in Appendix~\ref{app:derivative}.

Equation~\eqref{eq:totalerror} shows that the total error between the HRBE and truth solutions is bounded by the sum of the RB approximation and hyperreduction errors. Consequently, to ensure that the RB approximation remains the dominant source of error at the system level, it suffices to choose the hyperreduction accuracy such that $\wtilde{\varepsilon}_{\rb, \wc} \ll \varepsilon_{\rb, \wc}$ for all $\wc \in \wcset$.\footnote{We set $\wtilde{\varepsilon}_{\rb, \wc}$ to be much smaller than $\varepsilon_{\rb, \wc}$ to ensure that the hyperreduction is smaller than the RB approximation error, but not so excessively so that the hyperreduction is still efficient. In practice, the ratio of the two residuals is set to $\calO(10)$--$\calO(100)$.}

\subsubsection{Coupling hyperreduction and POD fidelities for component-wise training}
\label{subsubsection:coupling}
We now propose a strategy to couple hyperreduction and POD fidelities in the component-wise offline training. As outlined in Section~\ref{subsec:rbtraining}, each $\wc \in \wcset$ is associated with a family of RB spaces $\{\what{\mathcal{V}}_{\rb,\wc}^{\mathbf{f}}\}_{\mathbf{f} \in \mathcal{F}_\wc}$. Following~\eqref{eq:epsilonhc} and~\eqref{eq:epsilonrbc}, we introduce $\varepsilon_{\rb, \wc}^{\mathbf{f}}$ and $\wtilde{\varepsilon}_{\rb, \wc}^{\mathbf{f}}$ as the RB approximation and hyperreduction errors for the RB space $\what{\mathcal{V}}_{\rb,\wc}^{\mathbf{f}}$, respectively. 

We first describe the procedure for approximating $\varepsilon_{\rb, \wc}^{\mathbf{f}}$ $\forall \wc \in \wcset$ and $\forall \mathbf{f} \in \mathcal{F}_\wc$. For each archetype component $\wc$, we consider the truth snapshot set $U_{h,\wc}^\train$ generated using Algorithm~\ref{alg:emptr}. We define $P_{\what{\mathcal{V}}_{\rb,\wc}^{\mathbf{f}}}: \what{\mathcal{V}}_{h,\wc} \to \what{\mathcal{V}}_{\rb,\wc}^{\mathbf{f}}$ as the $H^1(\what{\Omega}_\wc)$ projection operator onto the RB space $\what{\mathcal{V}}_{\rb,\wc}^{\mathbf{f}}$. Then, using~\eqref{eq:epsilonhc}, we estimate the RB approximation error $\varepsilon_{\rb, \wc}^{\mathbf{f}}$ as the maximum dual norm of the component's truth residual over all projected snapshot solutions; i.e,
\aequation{\nonumber
  \varepsilon_{\rb, \wc}^{\mathbf{f}} \approx \max_{u_{h,\wc}^\train \in U_{h,\wc}^\train} \sup_{v \in \what{\mathcal{V}}_{h,\wc}} \frac{|\what{R}_{h,\wc}( P_{\what{\mathcal{V}}_{\rb,\wc}^{\mathbf{f}}} (u_{h,\wc}^\train), {v}; \mu_\wc^\train) |}{\| v \|_{\what{\mathcal{V}}_\wc}},
}
where $\mu_\wc^\train \in \Xi_\wc^\train$ is the parameter associated with the snapshot solution $u_{h,\wc}^\train$. 

Subsequently, following the error bound in~\eqref{eq:totalerror}, $\forall \wc \in \wcset$, we set $\wtilde{\varepsilon}_{\rb, \wc}^{\mathbf{f}} \ll \varepsilon_{\rb, \wc}^{\mathbf{f}}$ $\forall \mathbf{f} \in \mathcal{F}_\wc$ to ensure that the error between the HRBE and truth solutions is dominated by the RB approximation error, rather than by the hyperreduction error. Therefore, for each component $\wc$ and for each of its RB spaces $\what{\mathcal{V}}_{\rb,\wc}^{\mathbf{f}}$, we apply the component-wise EQP such that the constraint~\eqref{eq:epsilonrbc} is enforced with the computed hyperreduction error $\wtilde{\varepsilon}_{\rb, \wc}^{\mathbf{f}}$ following the procedure introduced in~\cite{ebrahimi2024hyperreduced}.

\begin{remark}\label{rem:epsilonrbc}
  We note that for~\eqref{eq:totalerror} to hold $\forall \mu \in \mathcal{D}$, the RB approximation and hyperreduction errors, $\varepsilon_{\rb, \wc}^{\mathbf{f}}$ and $\wtilde{\varepsilon}_{\rb, \wc}^{\mathbf{f}}$, must ensure that the constraints~\eqref{eq:epsilonhc} and~\eqref{eq:epsilonrbc} are satisfied over the entire parameter domain. However, this is not feasible, as the RB solution $u_\rb(\mu)$ for all $\mu \in \mathcal{D}$ are not available during the offline training. Consequently, we enforce these constraints only over a representative set of training snapshots. This approach assumes that the training set constructed using Algorithm~\ref{alg:emptr} sufficiently captures the space of solutions that may be encountered during the online phase.
\end{remark}
 
\subsection{Component-wise error contraction estimation}
\label{subsec:errorcontraction}
We now present a procedure to estimate the error contraction factors required in the adaptive refinement strategy. For Proposition~\ref{prop:errestimator} to hold, the error contraction factors must satisfy the property that, for \emph{any} system $\sys$ composed of the trained archetype components and for \emph{any} parameter value $\mu \in \mathcal{D}$, the error between the truth and HRBE solutions on each instantiated component $c \in \sys$ contracts by a factor of $\eta_{M(c),{\mathbf{f}}_c}$, $\mathbf{f}_c \in \bar{\mathcal{F}}_{M(c)}$, when its RB space is refined from ${\mathcal{V}}_{\rb,M(c)}^{\mathbf{f}_c}$ to the next higher-fidelity RB space ${\mathcal{V}}_{\rb,M(c)}^{\mathbf{f}_c'}$, $\mathbf{f}_c' = \mathbf{f}_c + \mathbf{1}$, during the adaptive refinement procedure. 

As discussed in Remark~\ref{rem:epsilonrbc}, it is not feasible to guarantee this contraction behavior for all systems and parameters, since the full set of truth and HRBE solutions is not available during the offline training. Consequently, we estimate the error contraction factors using a representative training set constructed using Algorithm~\ref{alg:emptr}. We assume that the training set sufficiently captures the range of solutions expected during the online phase and that the relative error reduction observed over the training data generalizes well to new system configurations and parameter instances. Therefore, we conservatively estimate the error contraction factors $\eta_{\wc,\mathbf{f}}$ for each archetype component $\wc \in \wcset$ and $\mathbf{f} \in \bar{\mathcal{F}}_\wc$ as the worst-case ratio observed over the training set; i.e., 
\aequation{\label{eq:etaetimate}
  \eta_{\wc,\mathbf{f}} \approx \max_{\mu \in \Xi_\wc^{\train}} \frac{\|u_{h,\wc}(\mu) - \wtilde{{u}}_{\rb,\wc}'(\mu) \|_{\what{\mathcal{V}}_\wc}}{\| u_{h,\wc}(\mu) - \wtilde{u}_{\rb,\wc}(\mu) \|_{\what{\mathcal{V}}_\wc}}.
}
Here, $u_{h,\wc}(\mu) \in \what{\mathcal{V}}_{h,\wc}$ is the truth solution obtained by solving~\eqref{eq:truthref}, and $\wtilde{{u}}_{\rb,\wc}(\mu) \in {\mathcal{V}}_{\rb,\wc}^{\mathbf{f}}$ and $\wtilde{{u}}_{\rb,\wc}'(\mu) \in {\mathcal{V}}_{\rb,\wc}^{\mathbf{f}'}$, where $\mathbf{f}' = \mathbf{f} + \mathbf{1}$, are the HRBE solutions obtained by solving~\eqref{eq:nlhrbepr} using ${\mathcal{V}}_{\rb,\wc}^{\mathbf{f}}$ and ${\mathcal{V}}_{\rb,\wc}^{\mathbf{f}'}$, respectively, for a given $\mu \in \Xi_\wc^{\train}$. For the truth and HRBE problems, the Dirichlet boundary conditions on the component ports are obtained from the $H^1(\what{\Omega}_\wc)$ projection of the training snapshot solutions onto the corresponding truth and RB spaces, respectively.

\section{A priori error analysis}
\label{sec:error}
We now present the theoretical details for the derivation of the error bound~\eqref{eq:totalerror}. First, we introduce the BRR theorem~\cite{caloz1997numerical, veroy2005certified} and its variant to formulate the RB approximation and hyperreduction error bounds, respectively.

\begin{lemma}[BRR theorem]\label{lemma:brr}
  Let $G: \mathcal{W} \rightarrow \mathcal{W}'$ be a $C^1$ mapping from a Banach space $\mathcal{W}$ to its dual space $\mathcal{W}'$. We introduce  $v \in \mathcal{W}$ such that $DG(v) \in \mathcal{L}(\mathcal{W}; \mathcal{W}')$ is an isomorphism, where $\mathcal{L}(\mathcal{W}; \mathcal{W}')$ is the space of linear mappings from $\mathcal{W}$ to $\mathcal{W}'$, equipped with the norm $\| L \|_{\mathcal{L}(\mathcal{W}; \mathcal{W}')} \equiv \sup_{z \in \mathcal{W}} \| Lz \|_{\mathcal{W}'} / \|z \|_{\mathcal{W}}$ for all $L \in \mathcal{L}(\mathcal{W}; \mathcal{W}')$. We further introduce
  \begin{align}
  \varepsilon &\equiv \left\|G(v)\right\|_{\mathcal{W}'}, \label{eq:varepsilon}\\
  \delta &\equiv \left\|DG^{-1}(v)\right\|_{\mathcal{L}(\mathcal{W}'; \mathcal{W})}, \label{eq:delta}\\
  T(\alpha) &\equiv \sup_{z \in \bar{B}(v, \alpha)}  \left\| DG(v)  - DG(z) \right\|_{\mathcal{L}(\mathcal{W}; \mathcal{W}')},\label{eq:Talpha}
  \end{align}
  where $\bar{B}(v, \alpha) \equiv \{ z \in \mathcal{W}: \left\| z - v \right\|_{\mathcal{W}} \leq \alpha \}$. Assume $2 \delta T(2 \delta \varepsilon) \leq 1$. Then, there exists a unique $w \in \mathcal{W}$ such that $G(w) = 0$ in the ball $\bar{B}(v, 2 \delta \varepsilon)$ and $DG(w) \in \mathcal{L}(\mathcal{W}; \mathcal{W}')$ is invertible and satisfies
    
  \begin{equation}\nonumber
  \left\| DG^{-1}(w) \right\|_{\mathcal{L}(\mathcal{W}'; \mathcal{W})} \leq 2 \delta.
  \end{equation}
  Additionally, $\forall z \in \bar{B}(v, 2 \delta \varepsilon)$, 
  \begin{equation}\label{eq:brrerr}
      \left\| z -w \right\|_{\mathcal{W}} \leq 2 \delta \left\| G(z) \right\|_{\mathcal{W}'}.
  \end{equation}
\end{lemma}
\begin{proof}
  See~\cite{caloz1997numerical}.
\end{proof}

\begin{lemma}\label{lemma:brreqp}
  Let $G: \mathcal{W} \rightarrow \mathcal{W}'$ and $\wtilde{G}: \mathcal{W} \rightarrow \mathcal{W}'$ be two $C^1$ mappings from a Banach space $\mathcal{W}$ to its dual space $\mathcal{W}'$. We introduce  $w \in \mathcal{W}$ such that $G(w) = 0$. We assume $DG(w) \in \mathcal{L}(\mathcal{W}; \mathcal{W}')$ is an isomorphism and define $\gamma \equiv \| DG^{-1}(w) \|_{\mathcal{L}(\mathcal{W}'; \mathcal{W})}$. We further introduce
  \begin{align}
  \wtilde{\varepsilon} &\equiv \|\wtilde{G}(w)\|_{\mathcal{W}'}, \label{eq:brreqpcond1}\\
  \wtilde{T}(\alpha) &\equiv \sup_{z \in \bar{B}(w, \alpha)}  \| DG(w) - D\wtilde{G}(z) \|_{\mathcal{L}(\mathcal{W}; \mathcal{W}')}. \label{eq:brreqpcond2}
  \end{align}
  We assume $2 \gamma \wtilde{T}(2 \gamma \wtilde{\varepsilon}) \leq 1$. Then, there exists a unique $\wtilde{w} \in \mathcal{W}$ such that $\wtilde{G}(\wtilde{w}) = 0$ in the ball $\bar{B}(w, 2 \gamma \wtilde{\varepsilon})$. Additionally, $\forall z \in \bar{B}(w, 2 \gamma \wtilde{\varepsilon})$ 
  \begin{equation}\label{eq:brreqperr}
      \| z - \wtilde{w} \|_{\mathcal{W}} \leq 2 \gamma \| \wtilde{G}(z) \|_{\mathcal{W}'}.
  \end{equation} 
\end{lemma}
\begin{proof}
  See Appendix~\ref{app:brreqpproof}.
\end{proof}

We now apply these lemmas to bound the RB approximation and hyperreduction errors. First, for all $v, w, z \in \mathcal{V}_{\rb}$ and $\mu \in \mathcal{D}$, we define $\wtilde{R}'_{\rb}(v, w, z; \mu)$ as the G\^{a}teaux derivative of the HRBE residual form $\wtilde{R}_{\rb}(\cdot,z;\mu)$ in~\eqref{eq:nlhrbepr} at $v$ in the direction of $w$. Explicit expression of $\wtilde{R}'_{\rb}(v,w,z; \mu)$ is provided in Appendix~\ref{app:derivative}. Next, we define the inf--sup constant
\aequation{\nonumber
    \wtilde{\beta}_\rb(v; \mu) &\equiv \inf_{w \in \mathcal{V}_\rb} \sup_{z \in \mathcal{V}_\rb} \frac{|\wtilde{R}'_\rb(v, w, z; \mu)|}{\| w\|_{\mathcal{V}} \| z\|_{\mathcal{V}}} \quadd && \forall v \in \mathcal{V}_\rb, \forall \mu \in \mathcal{D},
}
and present the following proposition.

\begin{proposition}[RB approximation error bound]\label{prop:rberror}
  For the system $\sys$ and given $\mu = (\mu_c)_{c \in \sys} \in \mathcal{D}$, we introduce
  \begin{align}
      T_h(\alpha) &\equiv \sup_{{z} \in \bar{B}({{u}_\rb(\mu)}, \alpha)} \sup_{v \in \mathcal{V}_h} \sup_{w \in \mathcal{V}_h} \frac{|  R'_h({u}_\rb(\mu), v, w; \mu) - R'_h(z, v, w; \mu)|}{\| v \|_{\mathcal{V}} \| w \|_{\mathcal{V}}}, \nonumber
  \end{align}
  where $u_\rb(\mu) \in \mathcal{V}_\rb$ is the system's port-reduced RB solution as in~\eqref{eq:nlrbpr}. We assume the condition~\eqref{eq:epsilonhc} hold, $\beta_h(u_\rb(\mu); \mu) > 0$, and
  \aequation{\label{eq:RBcond}
      T_h(\frac{2}{\beta_h(u_\rb(\mu); \mu)} \sum_{c \in \sys} \varepsilon_{\rb, M(c)}) &\leq \frac{\beta_h(u_\rb(\mu); \mu)}{2}.
  }
  
  Then, the truth solution ${u}_h(\mu) \in \mathcal{V}_h$ in~\eqref{eq:truthref} is in the ball $\bar{B}(u_\rb(\mu), 2 \sum_{c \in \sys} \varepsilon_{\rb, M(c)} / \beta_h(u_\rb(\mu); \mu))$. Moreover, 
  \aequation
  {\label{eq:RBerr}
  \|u_h(\mu) - u_\rb(\mu) \|_{\mathcal{V}} \leq \frac{2}{\beta_h(u_\rb(\mu); \mu)} \sum_{c \in \sys} \varepsilon_{\rb, M(c)}.
  }
\end{proposition}
\begin{proof}
  We employ Lemma~\ref{lemma:brr} with $\mathcal{W} \equiv \mathcal{V}_h$, $G(z) \equiv R_h(z, \cdot; \mu) \in \mathcal{V}_h'$ for all $z \in \mathcal{V}_h$, and $v \equiv u_\rb(\mu)$. This leads to $\delta \equiv 1 / \beta_h(u_\rb(\mu); \mu)$ in~\eqref{eq:delta} and $T(\alpha) \equiv T_h(\alpha)$ in~\eqref{eq:Talpha}. Additionally, we have
  \aequation{\label{eq:RBGW}
    \left\|G(v)\right\|_{\mathcal{W}'} &= \|{R}_{h}(u_\rb(\mu), \cdot; \mu) \|_{\mathcal{V}_h'} = \sup_{z \in \mathcal{V}_h} \frac{|{R}_{h}(u_\rb(\mu), z; \mu) |}{\|z \|_{\mathcal{V}}} \leq \sum_{c \in \sys} \varepsilon_{\rb, M(c)},
  }
  where the inequality follows from the triangular inequality and applying~\eqref{eq:epsilonhc}. Therefore, we set $\varepsilon \equiv \sum_{c \in \sys} \varepsilon_{\rb, M(c)}$ in~\eqref{eq:varepsilon}. Subsequently, according to Lemma~\ref{lemma:brr}, if the assumption~\eqref{eq:RBcond} holds, since ${u}_h(\mu) \in \mathcal{V}_h$ uniquely satisfies ${R}_h({u}_h(\mu), z; \mu) = 0$ $\forall z \in \mathcal{V}_h$, it is in the ball $\bar{B}(u_\rb(\mu), 2 \sum_{c \in \sys} \varepsilon_{\rb, M(c)} / \beta_h(u_\rb(\mu); \mu))$. We finally appeal to~\eqref{eq:brrerr} and~\eqref{eq:RBGW} to obtain~\eqref{eq:RBerr}.
\end{proof}

Proposition~\ref{prop:rberror} shows that when the RB solution lies in a neighborhood where the nonlinearity measure $T_h(\cdot)$ satisfies the condition~\eqref{eq:RBcond}, the BRR error bound~\eqref{eq:RBerr} for nonlinear problems resembles the classical error bound for linear problems in terms of the dual norm of the residual and the inf--sup constant, except for an additional factor of 2.

\begin{proposition}[Hyperreduction error bound]\label{prop:hrberror} 
  For the system $\sys$ and given $\mu = (\mu_c)_{c \in \sys} \in \mathcal{D}$, we introduce
  \aequation
  {\nonumber
      T_\rb(\alpha) &\equiv \sup_{{z} \in \bar{B}({{u}_\rb(\mu)}, \alpha)} \sup_{v \in \mathcal{V}_\rb} \sup_{w \in \mathcal{V}_\rb} \frac{|  R'_\rb({u}_\rb(\mu), v, w; \mu) - \wtilde{R}'_\rb(z, v, w; \mu)|}{\| v \|_{\mathcal{V}} \| w \|_{\mathcal{V}}}.
  }
  We assume the condition~\eqref{eq:epsilonrbc} holds, and
  \aequation{\label{eq:HRBEcond}
      T_\rb(\frac{2}{\beta_\rb({u}_\rb(\mu); \mu)} \sum_{c \in \sys} \wtilde{\varepsilon}_{\rb, M(c)}) &\leq \frac{\beta_\rb({u}_\rb(\mu); \mu)}{2}.
  }
  (Since $\forall \mu \in \mathcal{D}$ the port-reduced RB problem is well-posed, $\beta_\rb(u_\rb(\mu); \mu) > 0$.) Then, there exists a unique solution $\wtilde{u}_\rb(\mu) \in \mathcal{V}_\rb$ in the ball $\bar{B}({u}_\rb(\mu), 2 \sum_{c \in \sys} \wtilde{\varepsilon}_{\rb, M(c)} / \beta_\rb({u}_\rb(\mu); \mu))$ such that $\wtilde{R}_\rb(\wtilde{u}_\rb(\mu), z; \mu) = 0$ $\forall z \in \mathcal{V}_\rb$. Moreover, 
  \aequation
  {\label{eq:HRBEerr}
  \|\wtilde{u}_\rb(\mu) - {u}_\rb(\mu) \|_{\mathcal{V}} \leq \frac{2}{\beta_\rb({u}_\rb(\mu); \mu)} \sum_{c \in \sys} \wtilde{\varepsilon}_{\rb, M(c)}.
  }
\end{proposition}
\begin{proof}
  We apply Lemma~\ref{lemma:brreqp} with $\mathcal{W} \equiv \mathcal{V}_\rb$, $G(z) \equiv R_\rb(z, \cdot; \mu)$ and $\wtilde{G}(z) \equiv \wtilde{R}_\rb(z, \cdot; \mu)$ for all $z \in \mathcal{V}_\rb$, and $w \equiv u_\rb(\mu)$. This leads to $\gamma \equiv 1 / \beta_\rb({u}_\rb(\mu); \mu)$ and $\wtilde{T}(\alpha) \equiv T_\rb(\alpha)$ in~\eqref{eq:brreqpcond2}. Additionally, we have
  \aequation{\label{eq:RBGWtilde}
    \|\wtilde{G}(w)\|_{\mathcal{W}'} &= \|\wtilde{R}_{\rb}({u}_\rb(\mu), \cdot; \mu) \|_{\mathcal{V}_\rb'} = \sup_{z \in \mathcal{V}_\rb} \frac{|\wtilde{R}_{\rb}(u_\rb(\mu), z; \mu) |}{\|z \|_{\mathcal{V}}} \leq \sum_{c \in \sys} \wtilde{\varepsilon}_{\rb, M(c)},
  }
  where the inequality follows from the triangular inequality and the application of~\eqref{eq:epsilonrbc}. Hence, we set $\wtilde{\varepsilon} \equiv \sum_{c \in \sys} \wtilde{\varepsilon}_{\rb, M(c)}$ in~\eqref{eq:brreqpcond1}. Subsequently, if the assumption~\eqref{eq:HRBEcond} holds, due to Lemma~\ref{lemma:brreqp}, there exists a unique solution $\wtilde{u}_\rb(\mu) \in \mathcal{V}_\rb$ in the ball $\bar{B}({u}_\rb(\mu), 2 \sum_{c \in \sys} \wtilde{\varepsilon}_{\rb, M(c)} / \beta_\rb({u}_\rb(\mu); \mu))$ that satisfies $\wtilde{R}_\rb(\wtilde{u}_\rb(\mu), z; \mu) = 0$ $\forall z \in \mathcal{V}_\rb$. We finally appeal to~\eqref{eq:brreqperr} and~\eqref{eq:RBGWtilde} to obtain~\eqref{eq:HRBEerr}.
\end{proof}

We now apply these propositions to derive the total error bound~\eqref{eq:totalerror}.

\begin{corollary}[Total error bound] For the system $\sys$ and given $\mu = (\mu_c)_{c \in \sys} \in \mathcal{D}$, if the conditions of Propositions~\ref{prop:rberror} and \ref{prop:hrberror} hold and 
  \aequation{\label{eq:betamin}
    \beta(\mu) = \min (\beta_h({u}_\rb(\mu); \mu), \beta_\rb({u}_\rb(\mu); \mu)) \quadd \forall \mu \in \mathcal{D},
  }
then the total error between the truth and HRBE solutions satisfies the bound~\eqref{eq:totalerror}.
\end{corollary}
\begin{proof}
  We have
  \aequation{\nonumber
      \|u_h(\mu) - \wtilde{u}_\rb(\mu) \|_{\mathcal{V}} &\leq \|u_h(\mu) - u_\rb(\mu) \|_{\mathcal{V}} + \|\wtilde{u}_\rb(\mu) - u_\rb(\mu)\|_{\mathcal{V}} \\
      &\leq \frac{2}{\beta_h(u_\rb(\mu);\mu)} \sum_{c \in \sys} \varepsilon_{\rb, M(c)} + \frac{2}{\beta_\rb(u_\rb(\mu); \mu)} \sum_{c \in \sys} \wtilde{\varepsilon}_{\rb, M(c)} \\
      &\leq \frac{2}{\beta(\mu)} \sum_{c \in \sys} \left( \varepsilon_{\rb, M(c)} + \wtilde{\varepsilon}_{\rb, M(c)} \right),
  }
  where the first inequality follows from the triangle inequality, the second inequality follows from~\eqref{eq:RBerr} and~\eqref{eq:HRBEerr}, and the last inequality follows from~\eqref{eq:betamin}.
\end{proof}

\section{Case study: nonlinear thermal fin systems}
\label{sec:res}
We now present a case study to demonstrate the performance of the adaptive HRBE method. We consider nonlinear thermal fin systems composed of an aluminum alloy, with a nonlinear temperature-dependent thermal conductivity $k: [1, 300] \: \mathrm{K} \to [4.341, 177.868] \: \mathrm{W/K}$ \cite{nist_aluminum_3003}. The thermal conductivity is modeled as 
\begin{equation}\label{eq:conductivity}
  \log (k(v)) = \sum_{i=0}^7 k_i \: (\log (v))^i \quad \forall v \in [1, 300] \: \mathrm{K},
  \end{equation}
where the coefficients $k_i$, $i = 0,\dots, 7$, are listed in Table~\ref{table:coeffs}. The parameterized continuous residual form for the systems is expressed as
  \begin{equation}\nonumber
  R(w,v;\mu) =  \int_{\Omega(\mu)} \left(k(w) \nabla w \right) \cdot \nabla v \: dx - \int_{\Omega(\mu)} f(\mu) \: v \: dx \quad \forall w,v \in \mathcal{V},
  \end{equation} 
where $\mathcal{V} = \left\{ v \in H^1(\Omega(\mu)) : \ v|_{\Gamma_D} = 0 \right\}$, and $f: \mathcal{D} \to L^2(\Omega(\mu))$ is the volumetric source term, assumed to be constant within each component. The first integral in the residual form depends nonlinearly on the field variable $w$ and does not admit an affine decomposition.

\begin{table}
  \caption{\label{table:coeffs}Coefficients of the aluminum's thermal conductivity equation~\eqref{eq:conductivity}.}
  \begin{center}
  \begin{tabular}{lcccccccc}
  \toprule[0.1em]
  {Coefficient} & $k_0$ & $k_1$ & $k_2$ & $k_3$ & $k_4$ & $k_5$ & $k_6$ & $k_7$\\[0.2em]
  {Value (W/K)} & 0.637 & -1.144 & 7.462 & -12.691 & 11.917 & -6.187 & 1.639 & -0.173\\[0.2em]
  \bottomrule
  \end{tabular}
  \end{center}
\end{table}

\subsection{Archetype port and component library}
\label{subsec:library}
\begin{figure}
  \centering
  \includegraphics[width=1\textwidth]{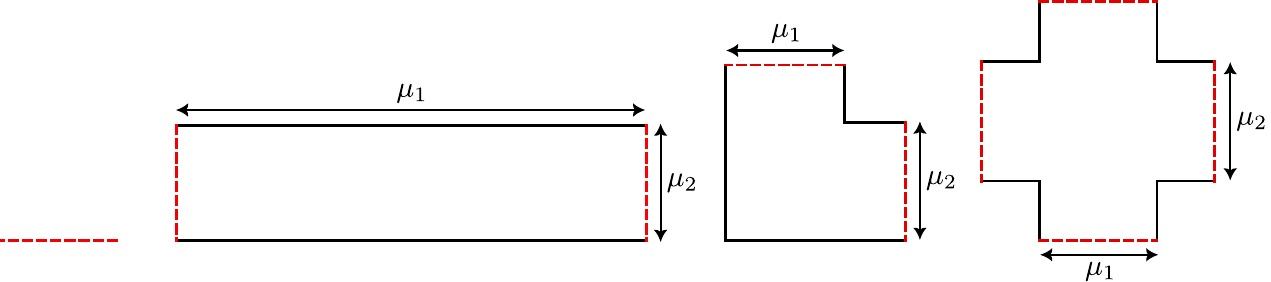}
          \caption{Archetype port and components in their reference domains. From left to right: line port, rod component, bracket component, and cross component. Ports are shown by red dashed lines.}\label{fig:lib}
  \end{figure}

Figure~\ref{fig:lib} shows the set of archetype components $\wcset = \{\mathrm{Rod}, \mathrm{Bracket}, \mathrm{Cross} \}$ considered in this study. The local ports of all archetype components are mapped from the same non-parameterized $17$-DoF archetype port depicted in Figure~\ref{fig:lib}. Hence, the set of archetype ports is $\wpset = \{ \mathrm{Line} \}$. The reference port has a length of $1$ cm and is discretized into eight quadratic line elements. The rod and bracket components each have two local ports, while the cross component has four. 

All archetype components are parameterized by two geometric parameters, $\mu_1$ and $\mu_2$. For all components, $\mu_1, \mu_2 \in [0.25, 1.5]$ cm, except for the rod component, where $\mu_1 \in [3, 6]$ cm. In their reference domains, $\mu_1 = \mu_2 = 1$ cm for all components, except for the rod component, where $\mu_1 = 4$ cm. Additionally, the components are parameterized by one physical parameter, $\mu_3 \in [0, 10]$ $\mathrm{W/cm^2}$, which characterizes the volumetric heat source. All components are discretized by quadratic triangular elements, which leads to $\mathcal{N}^{\bb}_{\mathrm{Rod}} = 691$, $\mathcal{N}^{\bb}_{\mathrm{Bracket}} = 703$, and $\mathcal{N}^{\bb}_{\mathrm{Cross}} = 1165$.

\subsubsection{Generation of snapshot solutions}

To generate the snapshot solutions for the archetype components and ports, we apply Algorithm~\ref{alg:emptr}. For each component, we create $N_{\sample} = 100$ sample subsystems by connecting it to other components in the library through its ports with a connection probability of $\nu = 0.8$. We assign uniformly random parameter values to the components in the subsystems and impose uniformly random constant Dirichlet boundary conditions to their boundary global ports, with values ranging from $1$~K to $250$~K. We solve the truth problem~\eqref{eq:truthref} for each assembled subsystem using Newton's method. The resulting truth snapshot solutions are stored in the full-component and bubble snapshot sets $U_{h, \wc}^{\train}$ and $U_{h, \wc}^{\train,\bb}$ $\forall \wc \in \wcset$, respectively. Since there is only one archetype port, all local port solutions are mapped to the reference domain of the line port and stored in its corresponding snapshot set $U_{h, \mathrm{Line}}^{\train}$.

\subsubsection{Construction of RB spaces}
We consider four fidelity levels to construct $N_{\mathrm{Line}}^{\mathrm{fdl}} = 4$ hierarchical RB spaces for the line port. To this end, we apply POD to $U_{h, \mathrm{Line}}^{\train}$ with tolerances $\delta_{\mathrm{Line}, \mathrm{pod}}^1 = 0.1$, $\delta_{\mathrm{Line}, \mathrm{pod}}^2 = 0.01$, $\delta_{\mathrm{Line}, \mathrm{pod}}^3 = 0.001$, and $\delta_{\mathrm{Line}, \mathrm{pod}}^4 = 0.0001$, where each POD tolerance is defined as the ratio between the square root of the sum of the discarded eigenvalues and the square root of the sum of all eigenvalues of the correlation matrix. This procedure yields RB spaces of dimension one, two, three, and five, respectively, corresponding to the four fidelity levels.

Similarly, for the bubble space of the archetype components, we consider four fidelity levels, i.e., $N_{\wc, \bb}^{\mathrm{fdl}} = 4$, for each $\wc \in \wcset$. To this end, we construct four hierarchical RB spaces for the bubble spaces by applying POD to $U_{h, \wc}^{\train,\bb}$ $\forall \wc \in \wcset$ with tolerances $\delta_{\wc, \mathrm{pod}}^1 = 0.1$, $\delta_{\wc, \mathrm{pod}}^2 = 0.01$, $\delta_{\wc, \mathrm{pod}}^3 = 0.001$, and $\delta_{\wc, \mathrm{pod}}^4 = 0.0001$. The resulting bubble space RB dimensions $N_\wc^\bb$ $\forall \wc \in \wcset$ at different fidelity levels are summarized in Table~\ref{table:archsummary}. The RB size ranges from 2 for the coarsest rod component to 22 for the finest cross component.

\begin{table}
  \caption{\label{table:archsummary} Summary of the offline training of the archetype components. }
  \begin{center}
  \begin{tabular}{lcccc}
  \toprule[0.1em]
  {Component} & {Rod} & {Bracket} & {Cross}\\[0.2em]
  \hline
  $\mathcal{N}^{\bb}_{\wc}$ & 691 & 703 & 1165 \\[0.2em]
  $Q_{\wc}$ & 1968 & 2016 & 3456 \\[0.0em]
  \hline \\[-1.em]
  $N^{\bb}_{\wc} \: (\delta_{\wc, \mathrm{pod}}^1 = 0.1)$ & 2 & 2 & 3 \\[0.2em]
  $N^{\bb}_{\wc} \: (\delta_{\wc, \mathrm{pod}}^2 = 0.01)$ & 3 & 3 & 8 \\[0.2em]
  $N^{\bb}_{\wc} \: (\delta_{\wc, \mathrm{pod}}^3 = 0.001)$ & 4 & 6 & 14 \\[0.2em]
  $N^{\bb}_{\wc} \: (\delta_{\wc, \mathrm{pod}}^4 = 0.0001)$ & 7 & 9 & 22 \\[0.2em]
  \hline \\[-1.em]
  $\wtilde{Q}_{\wc, \mathbf{f}}$ for $\mathbf{f} = (1, \dots, 1)$ & 7 & 32 & 56 \\[0.2em]
  $\wtilde{Q}_{\wc, \mathbf{f}}$ for $\mathbf{f} = (2, \dots, 2)$ & 15 & 42 & 131 \\[0.2em]
  $\wtilde{Q}_{\wc, \mathbf{f}}$ for $\mathbf{f} = (3, \dots, 3)$ & 32 & 89 & 233 \\[0.2em]
  $\wtilde{Q}_{\wc, \mathbf{f}}$ for $\mathbf{f} = (4, \dots, 4)$ & 44 & 133 & 296 \\[0.2em]
  \hline
  $\min_{\mathbf{f} \in \bar{\mathcal{F}}_\wc }\eta_{\wc, \mathbf{f}}$ & 0.207 & 0.235 & 0.214 \\[0.2em]
  $\max_{\mathbf{f} \in \bar{\mathcal{F}}_\wc }\eta_{\wc, \mathbf{f}}$ & 0.985 & 0.997 & 0.998 \\[0.2em]
  $\operatorname{median}_{\mathbf{f} \in \bar{\mathcal{F}}_\wc }\eta_{\wc, \mathbf{f}}$ & 0.423 & 0.383 & 0.352 \\[0.0em]
  \bottomrule
  \end{tabular}
  \end{center}
\end{table}

\subsubsection{Computation of RQ rules}
\begin{figure}
  \centering
      \begin{subfigure}{0.49\textwidth}
          \includegraphics[width=\textwidth]{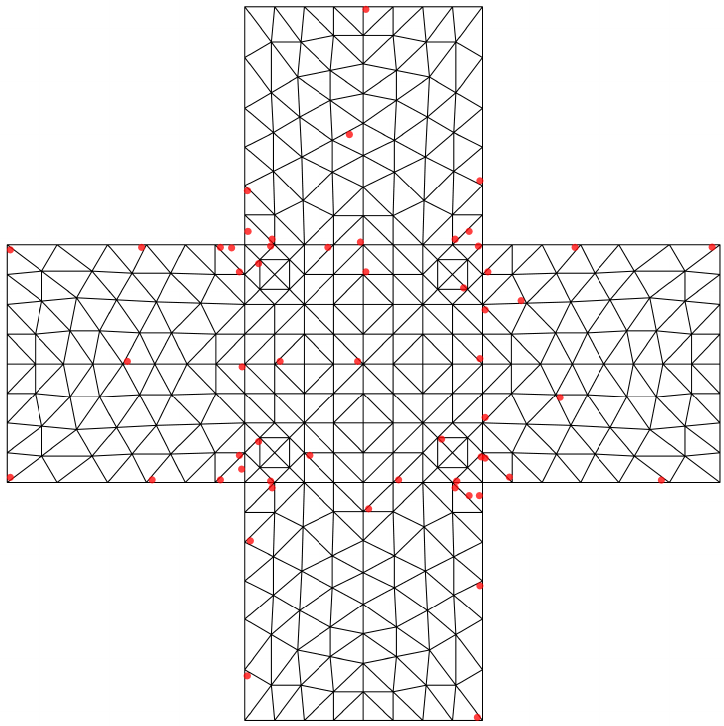}
          \caption{$\what{\mathcal{V}}_{\rb, \mathrm{Cross}}^{\mathbf{f}}$, $\mathbf{f} = (1, 1,1,1,1)$.}
      \end{subfigure}
      \begin{subfigure}{0.49\textwidth}
          \includegraphics[width=\textwidth]{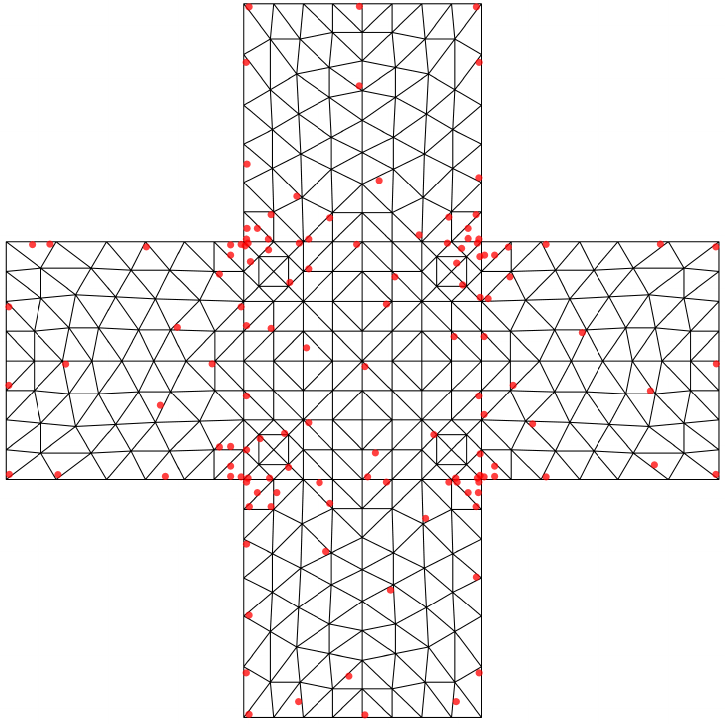}
          \caption{$\what{\mathcal{V}}_{\rb, \mathrm{Cross}}^{\mathbf{f}}$, $\mathbf{f} = (2, 2,2,2,2)$.}
      \end{subfigure}\\[0.5em]
      \begin{subfigure}{0.49\textwidth}
          \includegraphics[width=\textwidth]{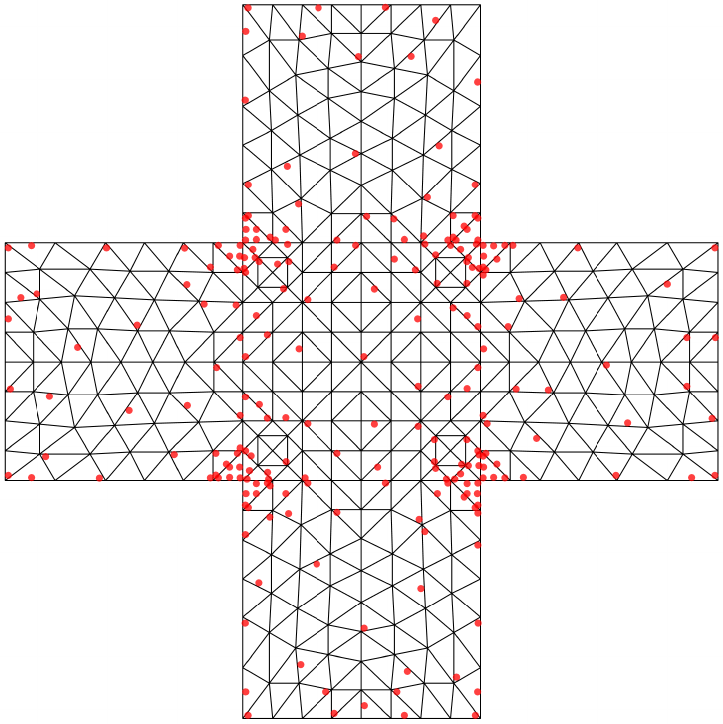}
          \caption{$\what{\mathcal{V}}_{\rb, \mathrm{Cross}}^{\mathbf{f}}$, $\mathbf{f} = (3, 3,3,3,3)$.}
      \end{subfigure}
      \begin{subfigure}{0.49\textwidth}
          \includegraphics[width=\textwidth]{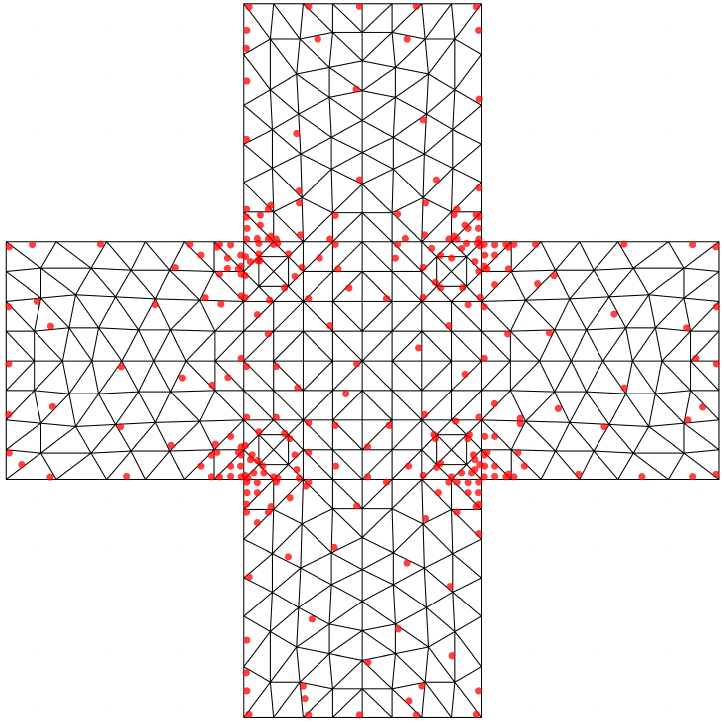}
          \caption{$\what{\mathcal{V}}_{\rb, \mathrm{Cross}}^{\mathbf{f}}$, $\mathbf{f} = (4, 4,4,4,4)$.}
      \end{subfigure}
      \caption{RQ points of the cross component for different multi-indexed RB spaces.}
      \label{fig:rqcross}
\end{figure}

We employ the component-wise EQP developed in~\cite{ebrahimi2024hyperreduced} to construct the RQ rules for the archetype components. To set the hyperreduction accuracies required for EQP, we first compute the RB approximation errors $\varepsilon_{\rb, \wc}^{\mathbf{f}}$ for each archetype component $\wc \in \wcset$ and each fidelity level $\mathbf{f} \in \mathcal{F}_\wc$, as described in Section~\ref{subsubsection:coupling}. To ensure that the hyperreduction error is smaller than the RB approximation error, at least for the training snapshot solutions, we set the hyperreduction accuracy $\wtilde{\varepsilon}_{\rb, \wc}^{\mathbf{f}}$ to $1\%$ of $\varepsilon_{\rb, \wc}^{\mathbf{f}}$ for each $\wc$ and ${\mathbf{f}}$. Figure~\ref{fig:rqcross} illustrates the resulting RQ rules for four selected multi-indexed RB spaces of the cross component. Table~\ref{table:archsummary} summarizes the number of RQ points for all $\wc \in \wcset$ across four selected multi-indexed RB spaces. As expected, the number of RQ points increases with the fidelity level of the RB spaces and complexity of the geometry. The number of RQ points ranges from 7 for the coarsest rod component to 296 for the finest cross component.

\begin{remark}
  For each archetype component $\wc \in \wcset$, we compute RQ rules for all of its $N_{\wc}^{\mathrm{fdl}}$ multi-indexed RB spaces. In our case study, this approach is feasible because the geometries are two-dimensional and the state variable (temperature) is scalar, which makes it computationally affordable to invoke the component-wise EQP many (i.e., $N_{\wc}^{\mathrm{fdl}}$) times. However, for problems involving higher-dimensional geometries and fields (e.g., three-dimensional hyperelasticity), the cost of performing EQP for all multi-indexed RB spaces can become prohibitive if $N_{\wc}^{\mathrm{fdl}}$ is large. In such cases, a practical strategy is to perform EQP only for a selected subset of the multi-indexed RB spaces. During adaptive refinement in the online phase, if an RQ rule for a given multi-indexed RB space is not available, we can employ the RQ rule associated with the closest higher-fidelity RB space for which EQP has been performed. This is justified because the set of constraint equations enforced in EQP for a smaller RB space is a subset of the set enforced for a larger RB space; thus, the RQ rule of a larger space remains valid for its subspaces. While this approach compromises some online computational efficiency---since more RQ points than strictly necessary are used for the smaller spaces---it reduces the offline training cost and the online memory footprint required to store the archetype library, and hence is a reasonable solution, especially when applying EQP for all multi-indexed spaces is computationally infeasible.
\end{remark}

\subsubsection{Estimation of error contraction factors}
We apply the procedure outlined in Section~\ref{subsec:errorcontraction} to conservatively estimate the error contraction factors for the archetype components. For each archetype component $\wc \in \wcset$ and each fidelity level $\mathbf{f} = (f_{\bb}, (f_{p})_{p \in \pset_{\wc}}) \in \bar{\mathcal{F}}_\wc$, we identify the next higher-fidelity RB space $\mathbf{f}' = \mathbf{f} + \mathbf{1}$, and estimate the error contraction factor $\eta_{\wc,\mathbf{f}}$ using~\eqref{eq:etaetimate}. Table~\ref{table:archsummary} reports the minimum, maximum, and median values of the error contraction factors for all archetype components. All contraction factors are strictly less than one, which confirms that refinement consistently improves the solution accuracy. Across all components, the smallest contraction factor, which indicates the most significant gain in accuracy, is observed when refining the RB space associated with $\mathbf{f} = (3, 1,\dots, 1)$ to $\mathbf{f}' = (4, 2,\dots, 2)$. In contrast, the largest contraction factor, which indicates the smallest improvement, is observed when refining the RB space that combines the highest-fidelity bubble space and all but one of the highest-fidelity port spaces to the fully enriched RB space, where all port spaces also reach their highest fidelity.

\begin{remark}
  Unlike the construction of RQ rules, which can be performed for a subset of multi-indexed RB spaces to manage offline computational cost, the estimation of error contraction factors must be carried out for all multi-indexed RB spaces. However, this requirement does not pose a significant burden, as the contraction factor estimation involves solving the truth and HRBE problems only at the component level. These are relatively inexpensive computations, which makes the overall cost of estimating contraction factors manageable.
\end{remark}

\subsection{Thermal fin systems}
\label{subsec:systems}
\begin{figure}
  \centering
   \begin{subfigure}{0.49\textwidth}
     \includegraphics[width=1.0\textwidth]{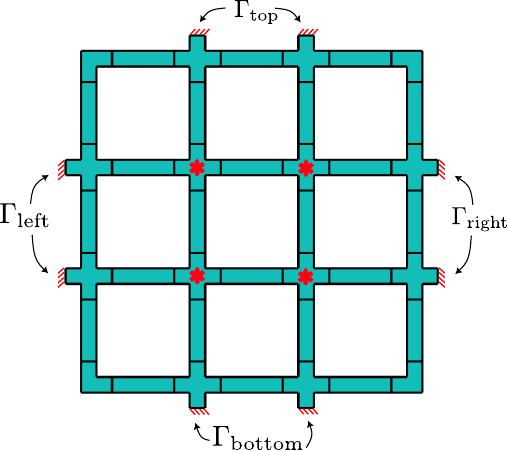}
     \caption{\label{fig:examplefins}}
   \end{subfigure}
   \begin{subfigure}{0.49\textwidth}
     \includegraphics[width=1.0\textwidth]{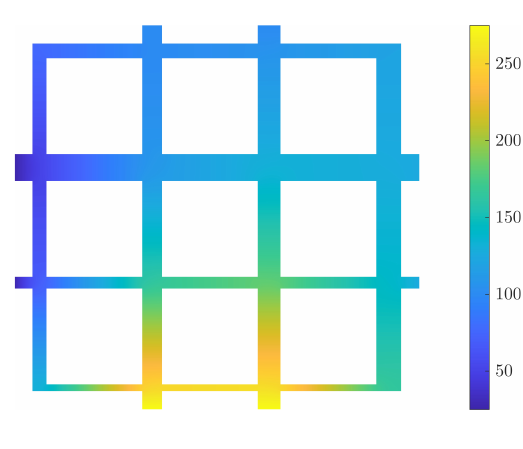}
     \caption{\label{fig:test1tr}}
   \end{subfigure}
   \caption{(a) An example fin system with $N_\fin = 3$. Components marked with red asterisks have a volumetric heat source. (b)~Truth temperature distribution for one instantiation of a system with $N_\fin = 3$.}
\end{figure}

We consider a family of thermal fin systems constructed from the archetype components described earlier. Figure~\ref{fig:examplefins} shows an example system consisting of 24 rod components, four bracket components, and 12 cross components. To simplify the parameterization, we focus on systems with an equal number of rod components along the horizontal and vertical directions. We characterize the topology of each system by a single parameter, $N_\fin$, which denotes the number of rod components per row/column in each direction. Thus, the system depicted in Figure~\ref{fig:examplefins} corresponds to a thermal fin system with $N_\fin = 3$. We study thermal fin systems with $N_\fin = 2$ to $8$. The number of instantiated components in each system is given by $N_\comp = (3 N_\fin + 1)(N_\fin + 1)$.

Table~\ref{table:truthvshrbe} reports the number of truth DoFs $\mathcal{N}_h$ and the number of truth quadrature points $Q_h = \sum_{c \in \sys} Q_{M(c)}$ for all values of $N_\fin$. The table also includes the number of HRBE DoFs ${N}_\rb^{(f)}$ and the corresponding number of RQ points $\wtilde{Q}_\rb^{(f)} = \sum_{c \in \sys} \wtilde{Q}_{M(c), (f, \dots, f)}$ for uniform fidelity level $f \in \{1, \dots, 4 \}$ across all components and ports. Even at the highest fidelity level, which is never actually reached in the case studies presented later, the HRBE problem involves $\approx 50\times$ fewer DoFs and $\approx 20\times$ fewer quadrature points than the truth model. This highlights the substantial potential for computational savings during the online phase.

\begin{table}
  \caption{\label{table:truthvshrbe} Summary of the number DoFs and quadrature points for truth and HRBE problems in the thermal fin systems.}
  \begin{center}
  \begin{tabular}{lccccccc}
  \toprule[0.1em]
  {$N_{\mathrm{fin}}$} & {$2$} & {$3$} & {$4$} & {$5$} & {$6$} & {$7$} & {$8$}\\[0.2em]
  \hline \\[-1.em]
  $\mathcal{N}_{h}$ & 17405 & 34328 & 56481 & 83864 & 116477 & 154320 & 197393 \\[0.2em]
  ${N}_\rb^{(1)}$ & 75 & 148 & 243 & 360 & 499 & 660 & 843 \\[0.2em]
  ${N}_\rb^{(2)}$ & 144 & 292 & 484 & 720 & 1000 & 1324 & 1692 \\[0.2em]
  ${N}_\rb^{(3)}$ & 226 & 456 & 754 & 1120 & 1554 & 2056 & 2626 \\[0.2em]
  ${N}_\rb^{(4)}$ & 370 & 748 & 1238 & 1840 & 2554 & 3380 & 4318 \\[0.2em]
  \hline \\[-1.em]
  $Q_{h}$ & 48960 & 96768 & 159360 & 236736 & 328896 & 435840 & 557568 \\[0.0em]
  $\wtilde{Q}_\rb^{(1)}$ & 492 & 968 & 1584 & 2340 & 3236 & 4272 & 5448 \\[0.2em]
  $\wtilde{Q}_\rb^{(2)}$ & 1003 & 2100 & 3519 & 5260 & 7323 & 9708 & 12415 \\[0.2em]
  $\wtilde{Q}_\rb^{(3)}$ & 1905 & 3920 & 6529 & 9732 & 13529 & 17920 & 22905 \\[0.2em]
  $\wtilde{Q}_\rb^{(4)}$ & 2540 & 5140 & 8508 & 12644 & 17548 & 23220 & 29660 \\[0.2em]
  \bottomrule
  \end{tabular}
  \end{center}
\end{table}

Four Dirichlet boundary conditions are imposed on the system boundaries as follows: (i) $u_{\mathrm{left}} = 25$~K on $\Gamma_{\mathrm{left}}$, (ii) $u_{\mathrm{right}} = 125$~K on $\Gamma_{\mathrm{right}}$, (iii) $u_{\mathrm{bottom}} = 275$~K on $\Gamma_{\mathrm{bottom}}$, and (iv) $u_{\mathrm{top}} = 100$~K on $\Gamma_{\mathrm{top}}$. Figure~\ref{fig:test1tr} illustrates the truth temperature distribution for one instantiation of a system with $N_{\fin} = 3$.

We parameterize each system by (i) a single length variable for the rod components, which is shared along the horizontal and vertical directions; (ii) $N_\fin + 1$ variables for component thicknesses along each direction; and (iii) $(N_\fin - 1)^2$ volumetric heat-source variables applied to the interior cross components. Therefore, each system is parameterized by a total of $N_\fin^2 + 4$ variables.

\subsection{Numerical results}
\subsubsection{First test: verification of the adaptation behavior for a localized heat source}
\begin{figure}
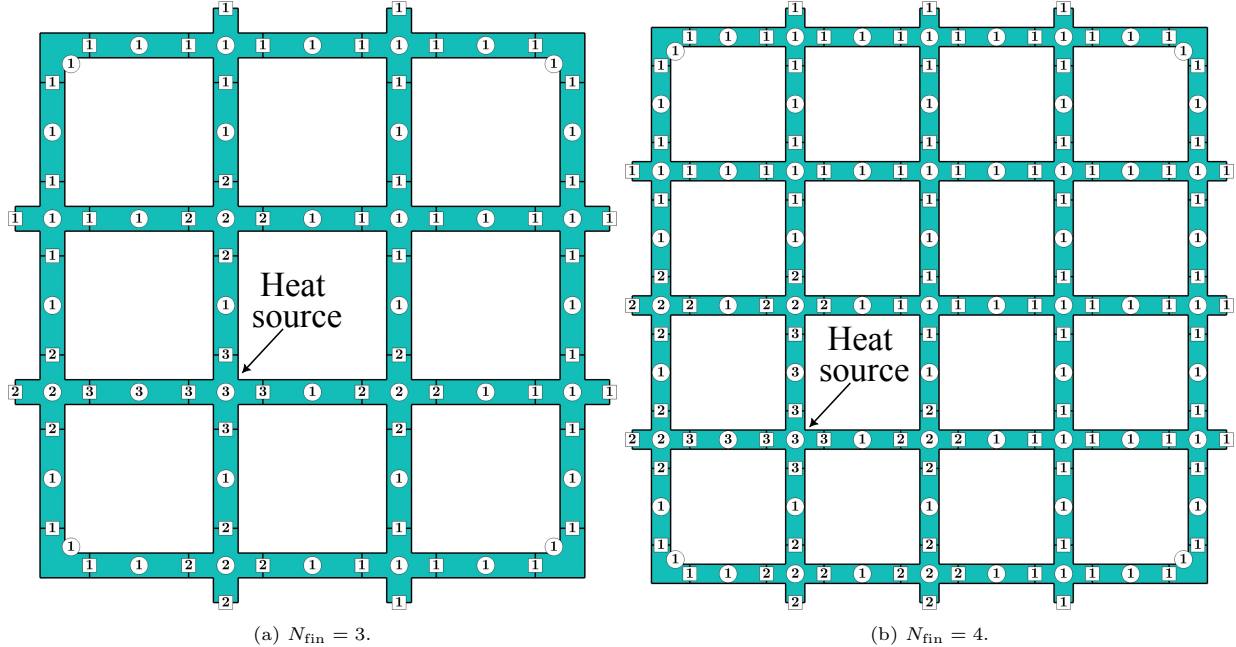

  \centering
   \begin{subfigure}[b]{0.49\textwidth}
     \includegraphics[width=1.0\textwidth]{3x3_7_fidelities_3.pdf}
     \caption{$N_\fin = 3$.}
   \end{subfigure}
   \begin{subfigure}[b]{0.49\textwidth}
     \includegraphics[width=1.0\textwidth]{4x4_7_fidelities_3.pdf}
     \caption{$N_\fin = 4$.}
   \end{subfigure}
   \caption{Final fidelity levels of the bubble and port RB spaces for systems with (a) $N_\fin = 3$ and (b) $N_\fin = 4$ for the first test with a localized heat source. Circled numbers indicate the fidelity levels of the bubble RB spaces, and numbers in squares indicate the fidelity levels of the port RB spaces.\label{fig:oneheat}}
\end{figure}

In this test, for all fin sizes, all geometric parameters are set to their reference values, and the volumetric heat sources are zero except for the bottom-left cross component, which is assigned a heat source of $10$~$\mathrm{W/cm^2}$. We use this problem with a localized heat source to demonstrate the adaptive refinement behavior. We apply Algorithm~\ref{alg:adaptiveref} with $N_{\mathrm{ref}} = 10$, $\Delta = 10$, and a relative $\mathcal{V}$-norm error tolerance of $1\%$. For all $N_\fin$ values, convergence is achieved after two adaptive refinement iterations. Figure~\ref{fig:oneheat} displays the final fidelity levels of the bubble and port RB spaces for systems with $N_\fin = 3$ and $N_\fin = 4$. Across all fin sizes, only the components located within two components of the heat source are enriched, which indicates that the adaptive refinement strategy effectively targets local enrichment and selectively refines components that have the greatest impact on solution accuracy.

\begin{figure}[t]
  \centering
  \includegraphics[width=0.49\textwidth]{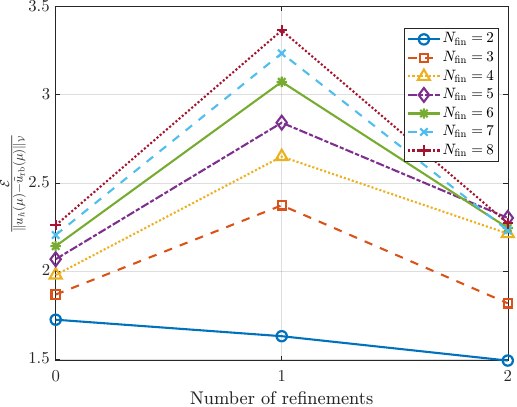}
  \caption{Effectivity of the error estimate~\eqref{eq:errestimator} for all fin sizes at different refinement iterations for the first test with a localized heat source.}\label{fig:eff_7}
\end{figure}

Figure~\ref{fig:eff_7} depicts the effectivity of the error estimate~\eqref{eq:errestimator} for all fin sizes at all adaptive refinement iterations. The effectivity is defined as the ratio between the error estimate and the actual error. In all cases, the effectivity is close to unity, which indicates that the proposed error estimate reliably approximates the true error and can be used to effectively guide the adaptive refinement process.

\subsubsection{Second test: comparison of the uniform and adaptive HRBE methods}
For this test, we compare the adaptive HRBE method against the uniform HRBE method, which refines the RB space fidelities uniformly across all components and ports in the system. This test uses the full set of parameters defined in Section~\ref{subsec:systems}. For each $N_\fin \in \{2,\dots, 8 \}$, we consider a test set $\Xi^{\test}_{N_\fin} \equiv \{ \mu^{\test}_{N_\fin, n} \in \mathcal{D} \}_{n=1}^5$, which consists of five parameter tuples sampled uniformly from their corresponding parameter space. We apply Algorithm~\ref{alg:adaptiveref} with $N_{\mathrm{ref}} = 10$, $\Delta = 20$, and a relative $\mathcal{V}$-norm error tolerance of $1\%$ to all cases.

\begin{figure}
  \centering
   \begin{subfigure}{0.49\textwidth}
     \includegraphics[width=1.0\textwidth]{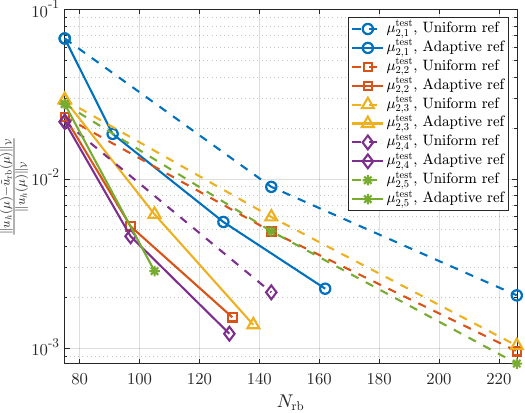}
     \caption{$N_\fin = 2$.}
   \end{subfigure}
   \begin{subfigure}{0.49\textwidth}
     \includegraphics[width=1.0\textwidth]{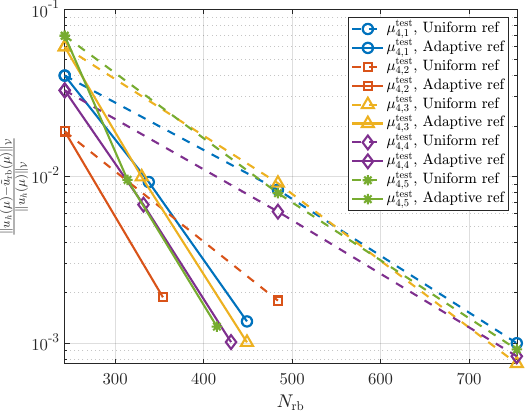}
     \caption{$N_\fin = 4$.}
   \end{subfigure}
   \begin{subfigure}{0.49\textwidth}
    \includegraphics[width=1.0\textwidth]{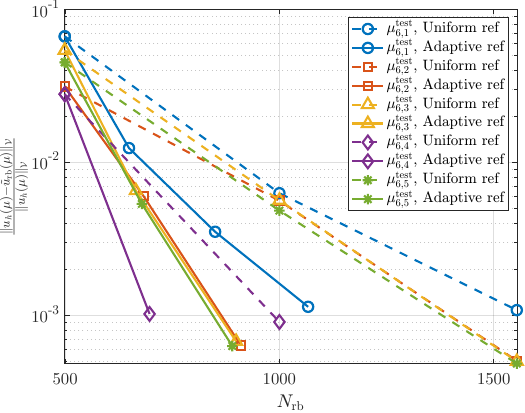}
    \caption{$N_\fin = 6$.}
  \end{subfigure}
  \begin{subfigure}{0.49\textwidth}
    \includegraphics[width=1.0\textwidth]{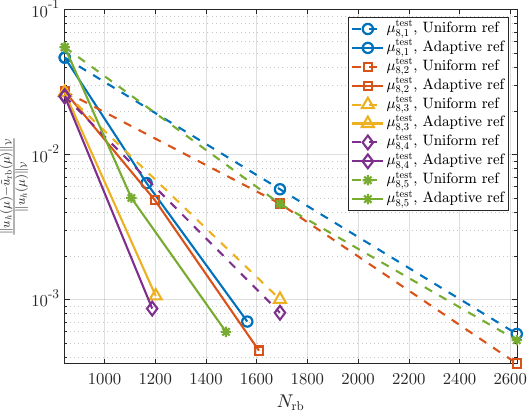}
    \caption{$N_\fin = 8$.}
  \end{subfigure}
   \caption{Convergence of the absolute error versus $N_\rb$ for the adaptive and uniform HRBE methods for $N_\fin \in \{2, 4, 6, 8\}$ and all test parameters for the second test.}\label{fig:convergence}
\end{figure}

Figure~\ref{fig:convergence} shows the evolution of the relative error between the truth and HRBE solutions versus $N_\rb$ for both the adaptive and uniform HRBE methods across $N_\fin \in \{2, 4, 6, 8\}$ and their test parameters.  For all fin sizes and test parameters, the adaptive HRBE method achieves the desired accuracy within two to four refinement iterations. In all cases, it reaches the target accuracy with fewer RB DoFs than the uniform HRBE method, which demonstrates that the adaptive strategy more efficiently directs computational resources toward the components that most influence solution accuracy. Specifically, for the $N_{\rm fin} = 8$ system, the adaptive HRBE method reduces the number of DoFs from $\mathcal{N}_h = 197393$ to $N_{\rm rb} = 1188$--$1609$ and the number of quadrature points from $Q_h = 557568$ to $\wtilde{Q}_\rb = \sum_{c \in \sys} \wtilde{Q}_{M(c), \mathbf{f}_c}= 8291$--$10463$, which corresponds to the reduction of $\approx 120$--$170\times$ and $\approx 50$--$70\times$, respectively.

\begin{figure}
  \centering
   \begin{subfigure}{0.48\textwidth}
     \includegraphics[width=1.0\textwidth]{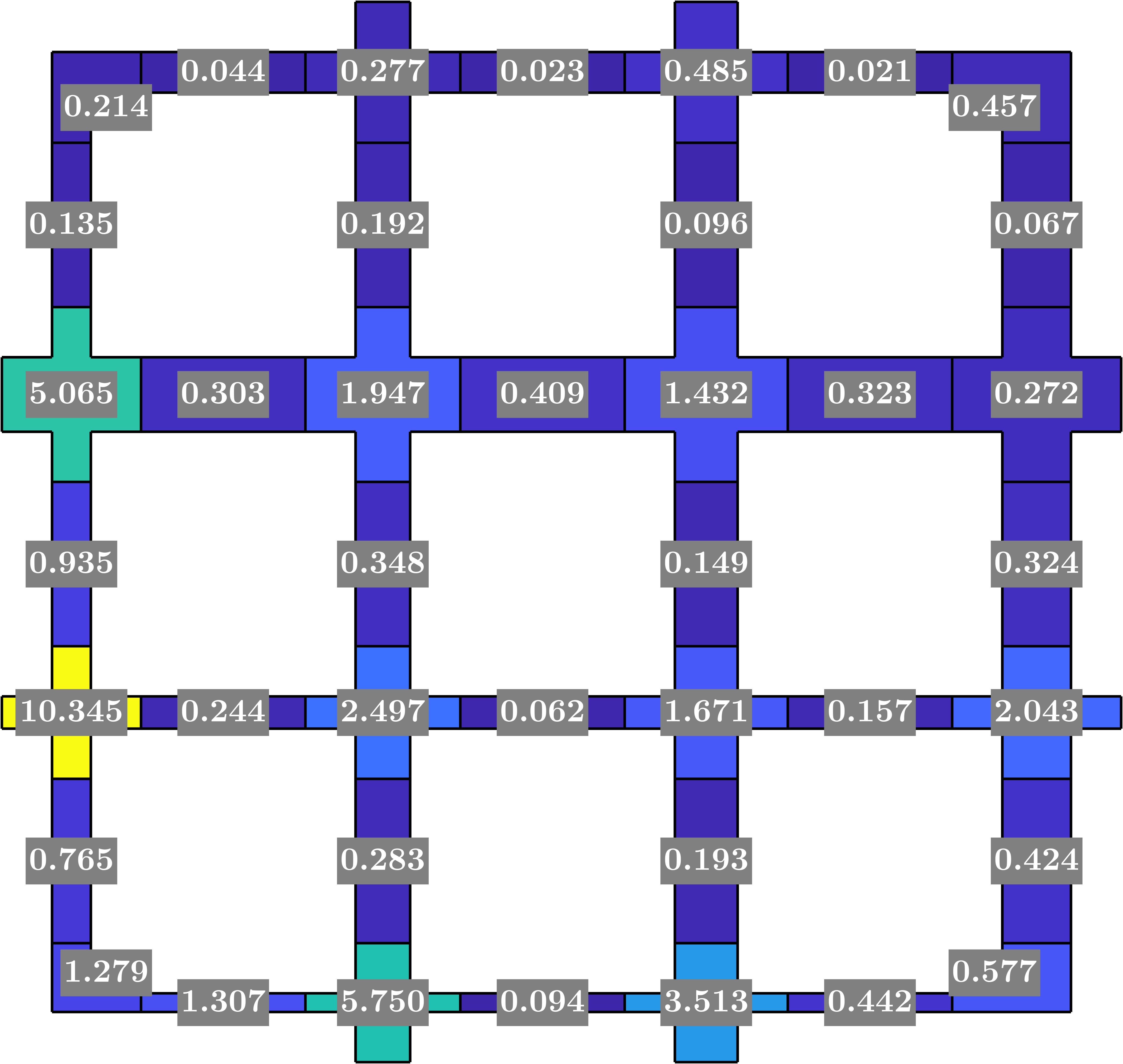}
     \caption{Component-wise error norms at the 1\textsuperscript{st} refinement iteration.}
   \end{subfigure}
   \begin{subfigure}{0.49\textwidth}
     \includegraphics[width=1.0\textwidth]{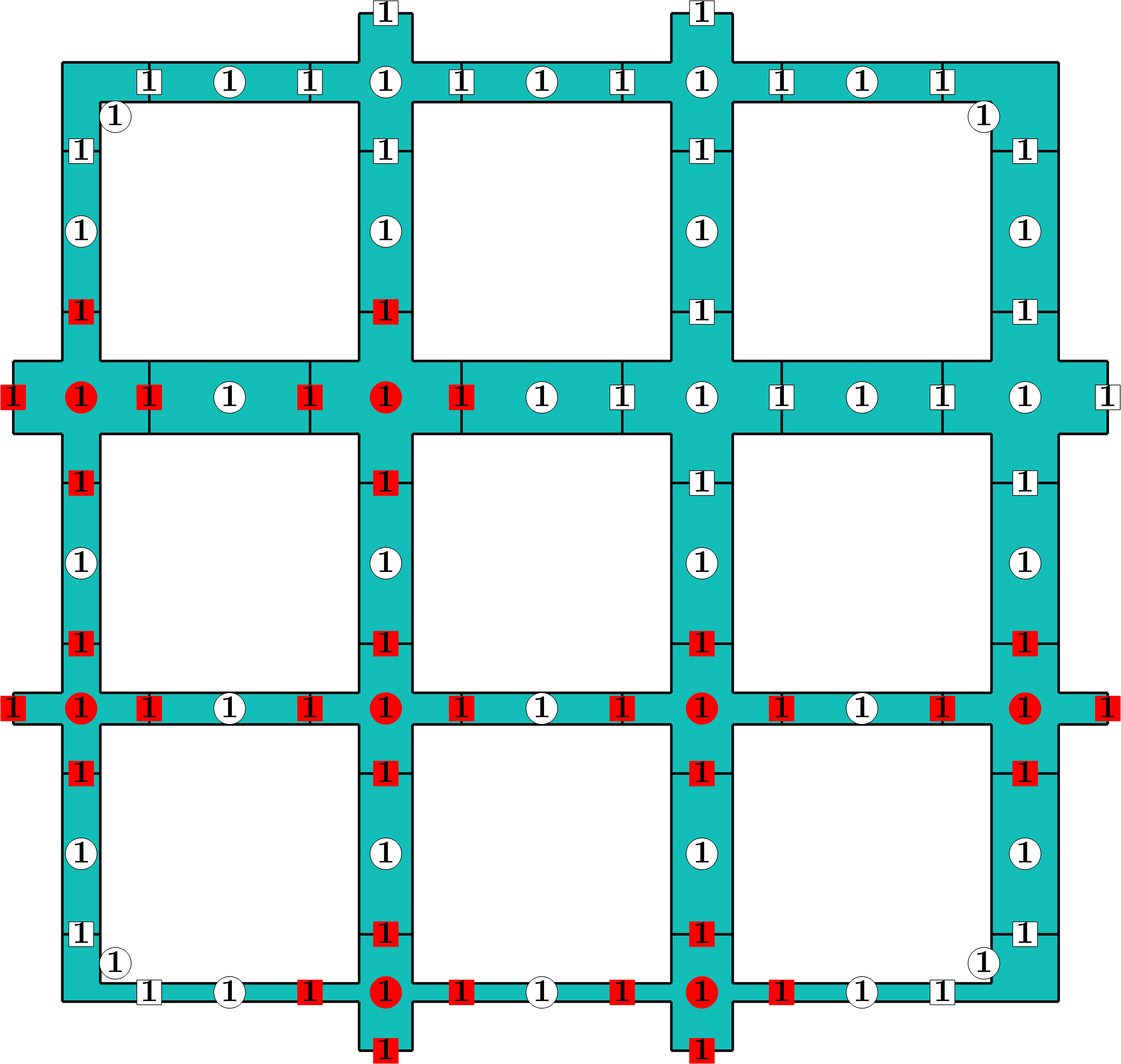}
     \caption{RB space fidelities at the 1\textsuperscript{st} refinement iteration.}
   \end{subfigure}
   \begin{subfigure}{0.48\textwidth}
    \includegraphics[width=1.0\textwidth]{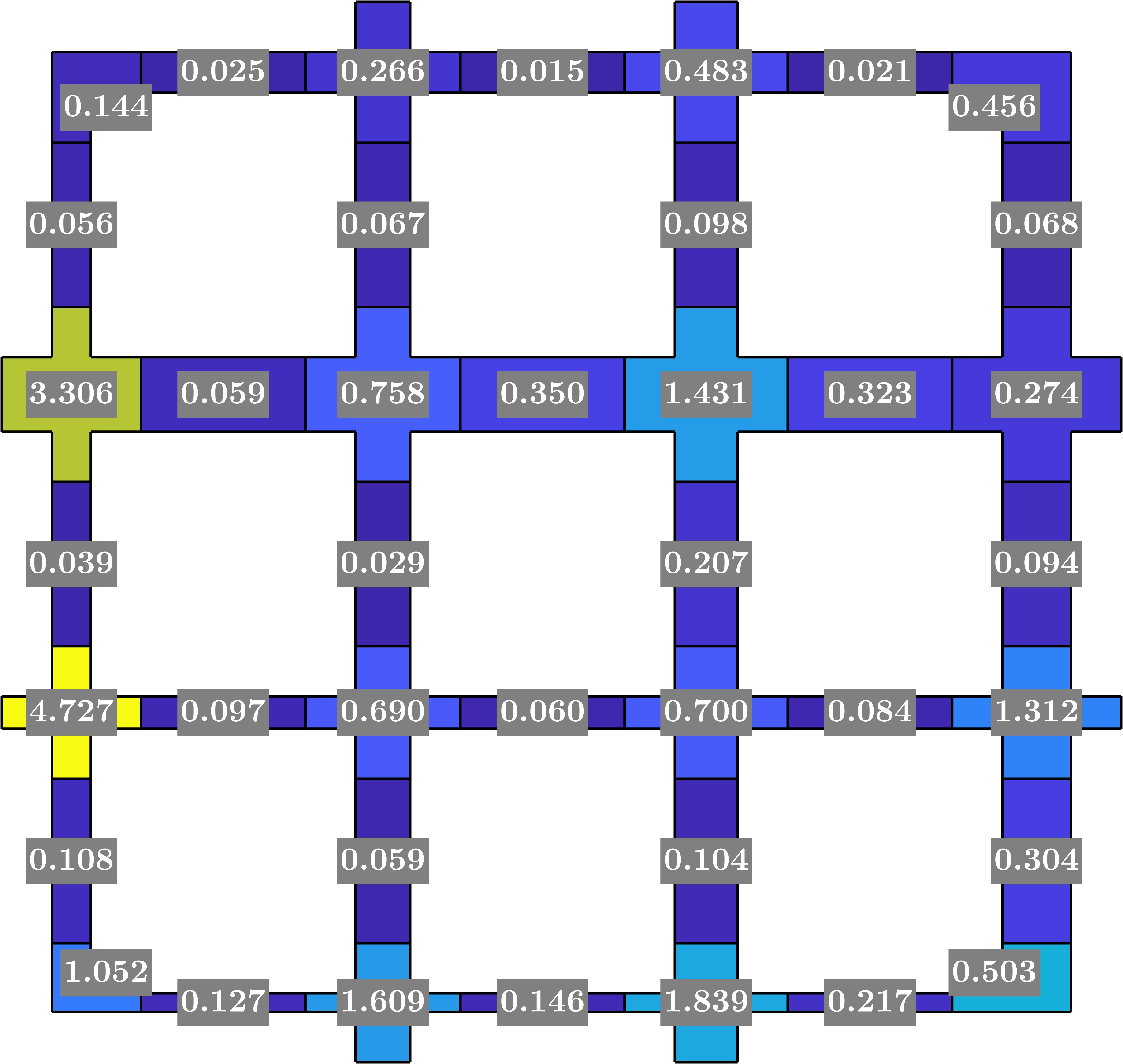}
    \caption{Component-wise error norms at the 2\textsuperscript{nd} refinement iteration.}
  \end{subfigure}
  \begin{subfigure}{0.49\textwidth}
    \includegraphics[width=1.0\textwidth]{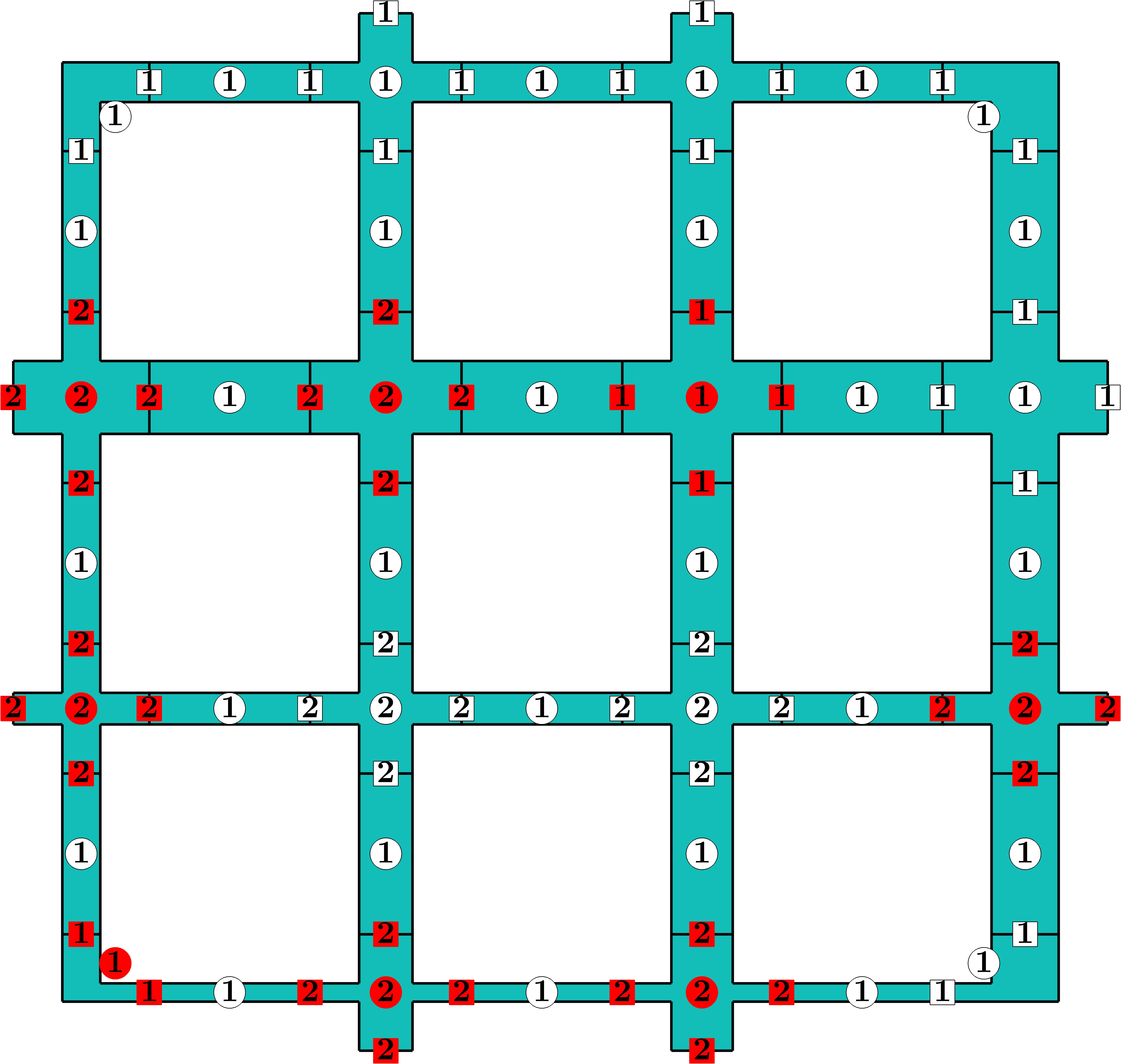}
    \caption{RB space fidelities at the 2\textsuperscript{nd} refinement iteration.}
  \end{subfigure}
   \caption{Evolution of component-wise errors and RB space fidelities at each adaptive refinement iteration for $N_\fin = 3$ and a selected test parameter for the second test. In the right column, red circles and squares mark the bubble and port RB spaces selected for refinement, respectively.}\label{fig:ref_3}
\end{figure}

Figure~\ref{fig:ref_3} depicts the evolution of the RB space fidelities and the actual component-wise error norms at each adaptive refinement iteration for $N_\fin = 3$ and a selected test parameter. While the figure presents a representative case, similar behavior is observed across all fin sizes and test parameters: in every case, the adaptive HRBE method identifies and refines the components with the largest errors, which confirms the strategy's local focus and effectiveness. Notably, the majority of components and ports remain at the lowest fidelity levels, which reduces computational cost without sacrificing accuracy.

\begin{figure}[t]
  \centering
  \includegraphics[width=0.65\textwidth]{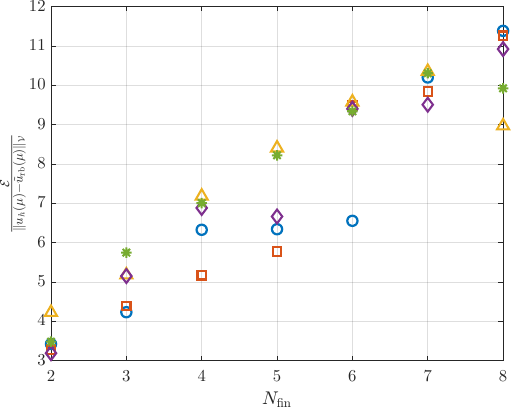}
  \caption{Effectivity of the error estimate~\eqref{eq:errestimator} for all fin sizes and test parameters for the second test.}\label{fig:eff_all}
\end{figure}

Figure~\ref{fig:eff_all} shows the effectivity of the error estimate~\eqref{eq:errestimator} for all fin sizes and test parameters. In all cases, the effectivity values range from 3 to 12, which confirms that the error estimate reliably tracks the true error across all system configurations and parameter samples. We observe a trend of increasing effectivity values with larger fin sizes. We suspect this is due to the use of conservative approximations of the error contraction factors, which in turn leads to increasingly conservative system-level error estimates as the number of components grows. Nonetheless, the variation in effectivity remains moderate, which indicates that the estimator retains its reliability even in larger systems.

\section{Conclusion}
\label{sec:conclusion}
In this work, we have developed an online-adaptive HRBE method for model reduction of parameterized, nonlinear systems composed of reusable components. The method constructs a library of archetype components during the offline phase, where each component is equipped with a family of hyperreduced models at multiple fidelity levels. During the online phase, these components are instantiated and assembled into a global system, where their fidelity levels are adaptively selected to satisfy a user-prescribed system-level error tolerance. This modular, reuse-driven strategy eliminates the need to retrain when the system topology changes and enables scalable reduced-order modeling of problems with high-dimensional parameter spaces.

A central contribution of this work is the development of an online adaptive refinement framework for components and ports informed by a hierarchical error estimator. The estimator compares solutions computed at successive fidelity levels and provides both local and system-level error estimates without requiring access to the truth model in the online phase. The local error indicators enable selective enrichment of the components that contribute most to the global error, which ensures that computational resources are used most effectively in the regions with the highest contribution to the error. 

Hyperreduction is performed using the component-wise EQP developed in~\cite{ebrahimi2024hyperreduced}. We apply the BRR theorem and its variant to develop a strategy that couples RB and hyperreduction fidelities during offline training to ensure that the hyperreduction error remains smaller than, but not excessively smaller than, the RB approximation error. This coupling reduces the complexity of online adaptation and guarantees that error control remains robust throughout the refinement process.

We demonstrate the proposed method on a family of nonlinear thermal fin systems. The results show that the adaptive HRBE method achieves $\calO(100)$ computational reduction while delivering accurate predictions with less than 1\% system-level error with respect to the truth model. Compared to uniform refinement, our adaptive approach more effectively targets high-error regions and achieves the desired accuracy with fewer DoFs. The system-level error estimator reliably tracks the actual error throughout the refinement process, which enables efficient and reliable online adaptation.

\newpage
\renewcommand\thesection{\Alph{section}}
\titleformat{\section}
{\normalfont\large\bfseries}
{Appendix~\thesection.}{1em}{}
\begin{appendices}
\section{Expressions of the G\^{a}teaux derivatives and residual forms}
\label{app:derivative}
We first provide the expressions for the G\^{a}teaux derivatives of the truth and HRBE residual forms. Using~\eqref{eq:truthref}, for all $v,w,z \in \mathcal{V}_h$ and $\mu \in \mathcal{D}$ we obtain
\begin{equation}\nonumber
    \begin{aligned}
    R_h'(v, w, z;\mu) = \sum_{c \in \sys} \sum_{q=1}^{Q_{M(c)}} \wrho_{M(c),q} \: \widehat{r}\:'_{M(c)} \Bigl( &\Bigl[ v^\bb_{c} + v^{\gamma}_{c} \Bigr] \circ \mathcal{G}_{c}(\cdot; \mu_c), \Bigl[ w^\bb_{c} + w^{\gamma}_{c} \Bigr] \circ  \mathcal{G}_{c}(\cdot; \mu_c), \\
    &\Bigl[ z^\bb_{c} + z^{\gamma}_{c} \Bigr] \circ \mathcal{G}_{c}(\cdot; \mu_c);\widehat{x}_{M(c),q},\mu_c \Bigr),
\end{aligned}
\end{equation}
where $\widehat{r}\:'_{\wc}(v,w,z;\widehat{x}_{\wc,q},\mu_\wc)$, $\wc \in \wcset$, is the G\^{a}teaux derivative of $\widehat{r}_{\wc}(\cdot,z;\widehat{x}_{\wc,q},\mu_c)$ at $v$ in the direction of~$w$. Similarly, using~\eqref{eq:nlhrbepr}, for all $v,w,z \in \mathcal{V}_\rb$ and $\mu \in \mathcal{D}$ we obtain
\begin{equation}\nonumber
  \begin{aligned}
  \wtilde{R}_\rb'(v, w, z;\mu) = \sum_{c \in \sys} \sum_{q=1}^{\wtilde{Q}_{M(c)}} \th{\rho}_{M(c),q} \: \widehat{r}\:'_{M(c)} \Bigl(&\Bigl[ v^\bb_{c} + v^{\gamma}_{c} \Bigr] \circ \mathcal{G}_{c}(\cdot; \mu_c), \Bigl[ w^\bb_{c} + w^{\gamma}_{c} \Bigr] \circ  \mathcal{G}_{c}(\cdot; \mu_c), \\
  &\hspace{0em} \Bigl[ z^\bb_{c} + z^{\gamma}_{c} \Bigr] \circ \mathcal{G}_{c}(\cdot; \mu_c);\th{x}_{M(c),q},\mu_c \Bigr),
\end{aligned}
\end{equation}
where $\widehat{r}\:'_{\wc}(v,w,z;\th{x}_{\wc,q},\mu_\wc)$, $\wc \in \wcset$, is the G\^{a}teaux derivative of $\widehat{r}_{\wc}(\cdot,z;\th{x}_{\wc,q},\mu_c)$ at $v$ in the direction of~$w$.

We now provide the expressions for the component residual forms. For all $v,w \in \what{\mathcal{V}}_{h,\wc}$ and $\mu \in \what{\mathcal{D}}_\wc$, the truth residual form for each archetype component $\wc \in \wcset$ is given by 
\begin{equation}\nonumber
  \begin{aligned}
  \what{R}_{h,\wc}(v, w;\mu) &= \sum_{q=1}^{Q_{\wc}} \wrho_{\wc,q} \widehat{r}_{\wc} \Big( v^\bb + v^{\gamma} ,w^\bb + w^{\gamma}; \widehat{x}_{\wc,q},\mu \Big).
\end{aligned}
\end{equation}
For all $v,w \in \what{\mathcal{V}}_{\rb,\wc}$ and $\mu \in \what{\mathcal{D}}_\wc$, the HRBE residual form for each archetype component $\wc \in \wcset$ is given by
\begin{equation}\nonumber
  \begin{aligned}
    \wtilde{\what{R}}_{\rb,\wc}(v, w;\mu) &= \sum_{q=1}^{\wtilde{Q}_{\wc}} \th{\rho}_{\wc,q} \widehat{r}_{\wc} \Big( v^\bb + v^{\gamma} ,w^\bb + w^{\gamma}; \th{x}_{\wc,q},\mu \Big).
\end{aligned}
\end{equation}

\section{Proof of Lemma~\ref{lemma:brreqp}}
\label{app:brreqpproof}
The proof proceeds similarly to that of the BRR theorem in~\cite{caloz1997numerical} with slight modifications. Define $H: \mathcal{W} \to \mathcal{W}$ as $H(z) = z - DG^{-1}(w) \wtilde{G}(z)$ $\forall z \in \mathcal{W}$. Accordingly, $z$ is a fixed point of $H(\cdot)$ if and only if $\wtilde{G}(z) = 0$. We have
\begin{equation}\nonumber
    \begin{aligned}
        H(z) - w &= z - w - DG^{-1}(w) \wtilde{G}(z) \\
        &= DG^{-1}(w)DG(w)(z - w) - DG^{-1}(w)(\wtilde{G}(z) - \wtilde{G}(w)) - DG^{-1}(w) \wtilde{G}(w) \\
        &= DG^{-1}(w)[DG(w)(z - w) - (\wtilde{G}(z) - \wtilde{G}(w))] - DG^{-1}(w) \wtilde{G}(w) \quadd \forall z \in \mathcal{W}.
    \end{aligned}
\end{equation}
We first show that $H(\cdot)$ maps $\bar{B}(w, 2 \gamma \wtilde{\varepsilon})$ into itself. The Taylor expansion of $\wtilde{G}(z)$ about $w$ is given by
\begin{equation}\nonumber
    \begin{aligned}
      \wtilde{G}(z) &= \wtilde{G}(w) + \int_{0}^1 D\wtilde{G}(w + t (z - w)) (z - w) dt.
    \end{aligned}
  \end{equation}
Thus, for any $z \in \bar{B}(w, 2 \gamma \wtilde{\varepsilon})$ we can write
\begin{equation}\nonumber
    \begin{aligned}
        \| H(z) - w \|_{\mathcal{W}} &= \| DG^{-1}(w)[DG(w)(z - w) - (\wtilde{G}(z) - \wtilde{G}(w))] - DG^{-1}(w) \wtilde{G}(w) \|_{\mathcal{W}} \\
        &\leq \Big\| DG^{-1}(w) \Big\|_{\mathcal{L}(\mathcal{W}'; \mathcal{W})} \: \Big\|\int_{0}^1 \Big(DG(w) - D\wtilde{G}(w + t (z - w)) \Big)  (w - z) dt + \wtilde{G}(w) \Big\|_{\mathcal{W}'}  \\
        &\leq \gamma (\wtilde{T}(2 \gamma \wtilde{\varepsilon}) 2 \gamma \wtilde{\varepsilon} + \wtilde{\varepsilon}) \leq 2 \gamma \wtilde{\varepsilon},
    \end{aligned}
\end{equation}
where the second inequality follows from~\eqref{eq:brreqpcond1}, \eqref{eq:brreqpcond2}, and the triangular inequality, and the last inequality follows from the assumption $2 \gamma \wtilde{T}(2 \gamma \wtilde{\varepsilon}) \leq 1$.

We now show that $H(\cdot)$ is a strict contraction mapping of $\bar{B}(w, 2 \gamma \wtilde{\varepsilon})$ into itself. For any $z \in \bar{B}(w, 2 \gamma \wtilde{\varepsilon})$, using the Taylor expansion of $\wtilde{G}(z)$ about $w$ we can write
\begin{equation}\nonumber
    \begin{aligned}
        \| H(z) - H(w) \|_{\mathcal{W}} &=  \Big\|  DG^{-1}(w) \int_{0}^1 \Big( DG(w) - D\wtilde{G}(w + t (z - w)) \Big)  (z - w) dt \Big\|_{\mathcal{W}} \\
        &\leq \gamma \wtilde{T}(2 \gamma \wtilde{\varepsilon}) \|z - w \|_{\mathcal{W}} \leq \frac{1}{2} \|z - w \|_{\mathcal{W}},
    \end{aligned}
\end{equation}
in which the first inequality follows from~\eqref{eq:brreqpcond2}, and the second inequality follows from the assumption $2 \gamma \wtilde{T}(2 \gamma \wtilde{\varepsilon}) \leq 1$. Hence, due to Banach's fixed-point theorem~\cite{khamsi2011introduction}, we conclude $H(\cdot)$ has a unique fixed point $\wtilde{w} \in \bar{B}(w, 2 \gamma \wtilde{\varepsilon})$. Consequently, $\wtilde{w} \in \mathcal{W}$ uniquely satisfies  $\wtilde{G}(\wtilde{w}) = 0$.

We now prove \eqref{eq:brreqperr}. For any $z \in \bar{B}(w, 2 \gamma \wtilde{\varepsilon})$ we can write
\begin{equation}\nonumber
  \begin{aligned}
      \wtilde{w} - z = \wtilde{w} - DG^{-1}(w) \wtilde{G}(\wtilde{w}) - z &= \wtilde{w} - z - DG^{-1}(w) \Big(\wtilde{G}(z) - \int_{0}^1 D\wtilde{G}(\wtilde{w} + t (z - \wtilde{w} )) (z - \wtilde{w}) dt \Big)\\
      &= DG^{-1}(w) (-\wtilde{G}(z) + \int_{0}^1 \Big(DG(w) - D\wtilde{G}(\wtilde{w} + t (z - \wtilde{w} )) \Big)  (\wtilde{w} - z) dt),
  \end{aligned}
\end{equation}
where the first equality follows from ``adding zero", and the second equality follows from the Taylor expansion of $\wtilde{G}(z)$ about $\wtilde{w}$. We then apply the triangular inequality and incorporate $2 \gamma \wtilde{T}(2 \gamma \wtilde{\varepsilon}) \leq 1$ to obtain~\eqref{eq:brreqperr}. \hspace*{\fill}\qedsymbol{}

\end{appendices}

\newpage
\renewcommand\refname{}
\section*{References}
\bibliographystyle{abbrv}
\bibliography{main}

\begin{thebibliography}{10}

\bibitem{barrault2004empirical}
M.~Barrault, Y.~Maday, N.~C. Nguyen, and A.~T. Patera.
\newblock An `empirical interpolation' method: application to efficient
  reduced-basis discretization of partial differential equations.
\newblock {\em Comptes Rendus Mathematique}, 339(9):667--672, 2004.

\bibitem{benner2015survey}
P.~Benner, S.~Gugercin, and K.~Willcox.
\newblock A survey of projection-based model reduction methods for parametric
  dynamical systems.
\newblock {\em SIAM Review}, 57(4):483--531, 2015.

\bibitem{bourquin1992component}
F.~Bourquin.
\newblock Component mode synthesis and eigenvalues of second order operators:
  discretization and algorithm.
\newblock {\em ESAIM: Mathematical Modelling and Numerical Analysis},
  26(3):385--423, 1992.

\bibitem{brezzi1980finite}
F.~Brezzi, J.~Rappaz, and P.-A. Raviart.
\newblock Finite dimensional approximation of nonlinear problems: {P}art~{I}:
  branches of nonsingular solutions.
\newblock {\em Numerische Mathematik}, 36(1):1--25, 1980.

\bibitem{buhr2018exponential}
A.~Buhr.
\newblock Exponential convergence of online enrichment in localized reduced
  basis methods.
\newblock {\em IFAC-PapersOnLine}, 51(2):302--306, 2018.

\bibitem{caloz1997numerical}
G.~Caloz and J.~Rappaz.
\newblock Numerical analysis for nonlinear and bifurcation problems.
\newblock {\em Handbook of Numerical Analysis}, 5:487--637, 1997.

\bibitem{chung2024train}
S.~W. Chung, Y.~Choi, P.~Roy, T.~Moore, T.~Roy, T.~Y. Lin, D.~T. Nguyen,
  C.~Hahn, E.~B. Duoss, and S.~E. Baker.
\newblock Train small, model big: {S}calable physics simulators via reduced
  order modeling and domain decomposition.
\newblock {\em Computer Methods in Applied Mechanics and Engineering},
  427:117041, 2024.

\bibitem{chung2024scaled}
S.~W. Chung, Y.~Choi, P.~Roy, T.~Roy, T.~Y. Lin, D.~T. Nguyen, C.~Hahn, E.~B.
  Duoss, and S.~E. Baker.
\newblock Scaled-up prediction of steady {N}avier-{S}tokes equation with
  component reduced order modeling.
\newblock {\em arXiv preprint arXiv:2410.21534}, 2024.

\bibitem{diaz2024fast}
A.~N. Diaz, Y.~Choi, and M.~Heinkenschloss.
\newblock A fast and accurate domain decomposition nonlinear manifold reduced
  order model.
\newblock {\em Computer Methods in Applied Mechanics and Engineering},
  425:116943, 2024.

\bibitem{ebrahimi2024hyperreduced}
M.~Ebrahimi and M.~Yano.
\newblock A hyperreduced reduced basis element method for reduced-order
  modeling of component-based nonlinear systems.
\newblock {\em Computer Methods in Applied Mechanics and Engineering},
  431:117254, 2024.

\bibitem{eftang2013port}
J.~L. Eftang and A.~T. Patera.
\newblock Port reduction in parametrized component static condensation:
  approximation and a posteriori error estimation.
\newblock {\em International Journal for Numerical Methods in Engineering},
  96(5):269--302, 2013.

\bibitem{eftang2014port}
J.~L. Eftang and A.~T. Patera.
\newblock A port-reduced static condensation reduced basis element method for
  large component-synthesized structures: approximation and a posteriori error
  estimation.
\newblock {\em Advanced Modeling and Simulation in Engineering Sciences},
  1(1):1--49, 2014.

\bibitem{grepl2007efficient}
M.~A. Grepl, Y.~Maday, N.~C. Nguyen, and A.~T. Patera.
\newblock Efficient reduced-basis treatment of nonaffine and nonlinear partial
  differential equations.
\newblock {\em ESAIM: Mathematical Modelling and Numerical Analysis},
  41(3):575--605, 2007.

\bibitem{hain2019hierarchical}
S.~Hain, M.~Ohlberger, M.~Radic, and K.~Urban.
\newblock A hierarchical a posteriori error estimator for the reduced basis
  method.
\newblock {\em Advances in Computational Mathematics}, 45(5):2191--2214, 2019.

\bibitem{hesthaven2016certified}
J.~S. Hesthaven, G.~Rozza, and B.~Stamm.
\newblock {\em Certified reduced basis methods for parametrized partial
  differential equations}.
\newblock Springer, 2016.

\bibitem{hoang2021domain}
C.~Hoang, Y.~Choi, and K.~Carlberg.
\newblock Domain-decomposition least-squares {Petrov}--{Galerkin} ({DD-LSPG})
  nonlinear model reduction.
\newblock {\em Computer Methods in Applied Mechanics and Engineering},
  384:113997, 2021.

\bibitem{huynh2013static}
D.~B.~P. Huynh, D.~J. Knezevic, and A.~T. Patera.
\newblock A static condensation reduced basis element method: approximation and
  a posteriori error estimation.
\newblock {\em ESAIM: Mathematical Modelling and Numerical Analysis},
  47(1):213--251, 2013.

\bibitem{huynh2013staticc}
D.~B.~P. Huynh, D.~J. Knezevic, and A.~T. Patera.
\newblock A static condensation reduced basis element method: {C}omplex
  problems.
\newblock {\em Computer Methods in Applied Mechanics and Engineering},
  259:197--216, 2013.

\bibitem{iollo2023one}
A.~Iollo, G.~Sambataro, and T.~Taddei.
\newblock A one-shot overlapping {Schwarz} method for component-based model
  reduction: application to nonlinear elasticity.
\newblock {\em Computer Methods in Applied Mechanics and Engineering},
  404:115786, 2023.

\bibitem{khamsi2011introduction}
M.~A. Khamsi and W.~A. Kirk.
\newblock {\em An introduction to metric spaces and fixed point theory}.
\newblock John Wiley \& Sons, 2011.

\bibitem{lovgren2006reduced}
A.~E. L{\o}vgren, Y.~Maday, and E.~M. R{\o}nquist.
\newblock A reduced basis element method for the steady {Stokes} problem.
\newblock {\em ESAIM: Mathematical Modelling and Numerical Analysis},
  40(3):529--552, 2006.

\bibitem{maday2002reduced}
Y.~Maday and E.~M. R{\o}nquist.
\newblock A reduced-basis element method.
\newblock {\em Journal of Scientific Computing}, 17(1):447--459, 2002.

\bibitem{maday2004reduced}
Y.~Maday and E.~M. Ronquist.
\newblock The reduced basis element method: application to a thermal fin
  problem.
\newblock {\em SIAM Journal on Scientific Computing}, 26(1):240--258, 2004.

\bibitem{nist_aluminum_3003}
{National Institute of Standards and Technology}.
\newblock Aluminum 3003-{F} ({UNS A93003}).
\newblock \url{https://www.nist.gov/mml/acmd/aluminum-3003-f-unsa93003}.
\newblock Accessed: August 2025.

\bibitem{ohlberger2015error}
M.~Ohlberger and F.~Schindler.
\newblock Error control for the localized reduced basis multiscale method with
  adaptive on-line enrichment.
\newblock {\em SIAM Journal on Scientific Computing}, 37(6):A2865--A2895, 2015.

\bibitem{Patera_2017_EQP}
A.~T. Patera and M.~Yano.
\newblock An {LP} empirical quadrature procedure for parametrized functions.
\newblock {\em Comptes Rendus Mathematique}, 355(11):1161--1167, 2017.

\bibitem{quarteroni2015reduced}
A.~Quarteroni, A.~Manzoni, and F.~Negri.
\newblock {\em Reduced basis methods for partial differential equations: an
  introduction}, volume~92.
\newblock Springer, 2015.

\bibitem{rozza2008reduced}
G.~Rozza, D.~B.~P. Huynh, and A.~T. Patera.
\newblock Reduced basis approximation and a posteriori error estimation for
  affinely parametrized elliptic coercive partial differential equations:
  application to transport and continuum mechanics.
\newblock {\em Archives of Computational Methods in Engineering}, 15(3):229,
  2008.

\bibitem{smetana2015new}
K.~Smetana.
\newblock A new certification framework for the port reduced static
  condensation reduced basis element method.
\newblock {\em Computer Methods in Applied Mechanics and Engineering},
  283:352--383, 2015.

\bibitem{smetana2016optimal}
K.~Smetana and A.~T. Patera.
\newblock Optimal local approximation spaces for component-based static
  condensation procedures.
\newblock {\em SIAM Journal on Scientific Computing}, 38(5):A3318--A3356, 2016.

\bibitem{smetana2023localized}
K.~Smetana and T.~Taddei.
\newblock Localized model reduction for nonlinear elliptic partial differential
  equations: localized training, partition of unity, and adaptive enrichment.
\newblock {\em SIAM Journal on Scientific Computing}, 45(3):A1300--A1331, 2023.

\bibitem{veroy2005certified}
K.~Veroy and A.~T. Patera.
\newblock Certified real-time solution of the parametrized steady
  incompressible {Navier--Stokes} equations: rigorous reduced-basis a
  posteriori error bounds.
\newblock {\em International Journal for Numerical Methods in Fluids},
  47(8-9):773--788, 2005.

\bibitem{wilson1974static}
E.~L. Wilson.
\newblock The static condensation algorithm.
\newblock {\em International Journal for Numerical Methods in Engineering},
  8(1):198--203, 1974.

\bibitem{yano2019lp}
M.~Yano and A.~T. Patera.
\newblock An {LP} empirical quadrature procedure for reduced basis treatment of
  parametrized nonlinear {PDEs}.
\newblock {\em Computer Methods in Applied Mechanics and Engineering},
  344:1104--1123, 2019.

\end{thebibliography}
\end{document}